\documentclass[a4paper,12pt]{article}
\usepackage{amsmath,amssymb,color}
  \usepackage{paralist}
  \usepackage{graphics} 
  \usepackage{epsfig} 
\usepackage{graphicx}  

\usepackage{lipsum}
\usepackage{amsfonts}
\usepackage{graphicx}
\usepackage{epstopdf}
\ifpdf
  \DeclareGraphicsExtensions{.eps,.pdf,.png,.jpg}
\else
  \DeclareGraphicsExtensions{.eps}
\fi

\title{
A Hybrid control approach to the route planning problem  for sailing boats
}

\author{Roberto Ferretti \footnote{Dipartimento di Matematica e Fisica, Universit\`a Roma Tre, L.go S. Leonardo Murialdo, 1, 00146 Roma, Italy. {\tt\small ferretti@mat.uniroma3.it}} \and Adriano Festa \footnote{INSA Rouen,
LMI laboratory, Avenue de l'Universit\'e, 76800 Saint-\'Etienne-du-Rouvray, France. {\tt\small adriano.festa@insa-rouen.fr} }}

\DeclareMathOperator{\tr}{tr}

\DeclareMathOperator{\R}{\mathbb{R}}
\DeclareMathOperator{\PP}{\mathbb{P}}
\DeclareMathOperator{\Q}{\mathbb{Q}}

\DeclareMathOperator{\N}{\mathbb{N}}

\newcommand{\E}{{\mathbb E}}
\newcommand{\I}{{\mathcal I}}

\newcommand{\RR}{{\mathbb R}}

\newcommand{\cN}{{\mathcal N}}

\newcommand{\cU}{{\mathcal U}}

\newcommand{\cT}{{\mathcal T}}

\newcommand{\cQ}{{\mathcal Q}}

\newtheorem{theorem}{\textbf{Theorem}}

\renewcommand{\phi}{\varphi}

\begin{document}

\maketitle

\begin{abstract}
We present an optimal hybrid control approach to the problem of stochastic route planning for sailing boats, especially in short course fleet races, in which minimum average time is an effective performance index. We show that the hybrid setting is a natural way of taking into account tacking/gybing maneuvers and other discrete control actions, and provide examples of increasing complexity to model the problem. Moreover, we carry out a numerical validation of the approach and show that results are in good agreement with theoretical and practical knowledge.
\end{abstract}


\section{Introduction}

In the last decades, the sport of sailing has experienced an increasing impact of new technologies, and notably of scientific computing. Among all computational problems relevant for sailing, we are interested here in \emph{route planning} and \emph{race strategy}, i.e., the optimization (with respect to a given performance index) of the yacht route.

\begin{figure}
\begin{center}
\includegraphics[height=5.4cm]{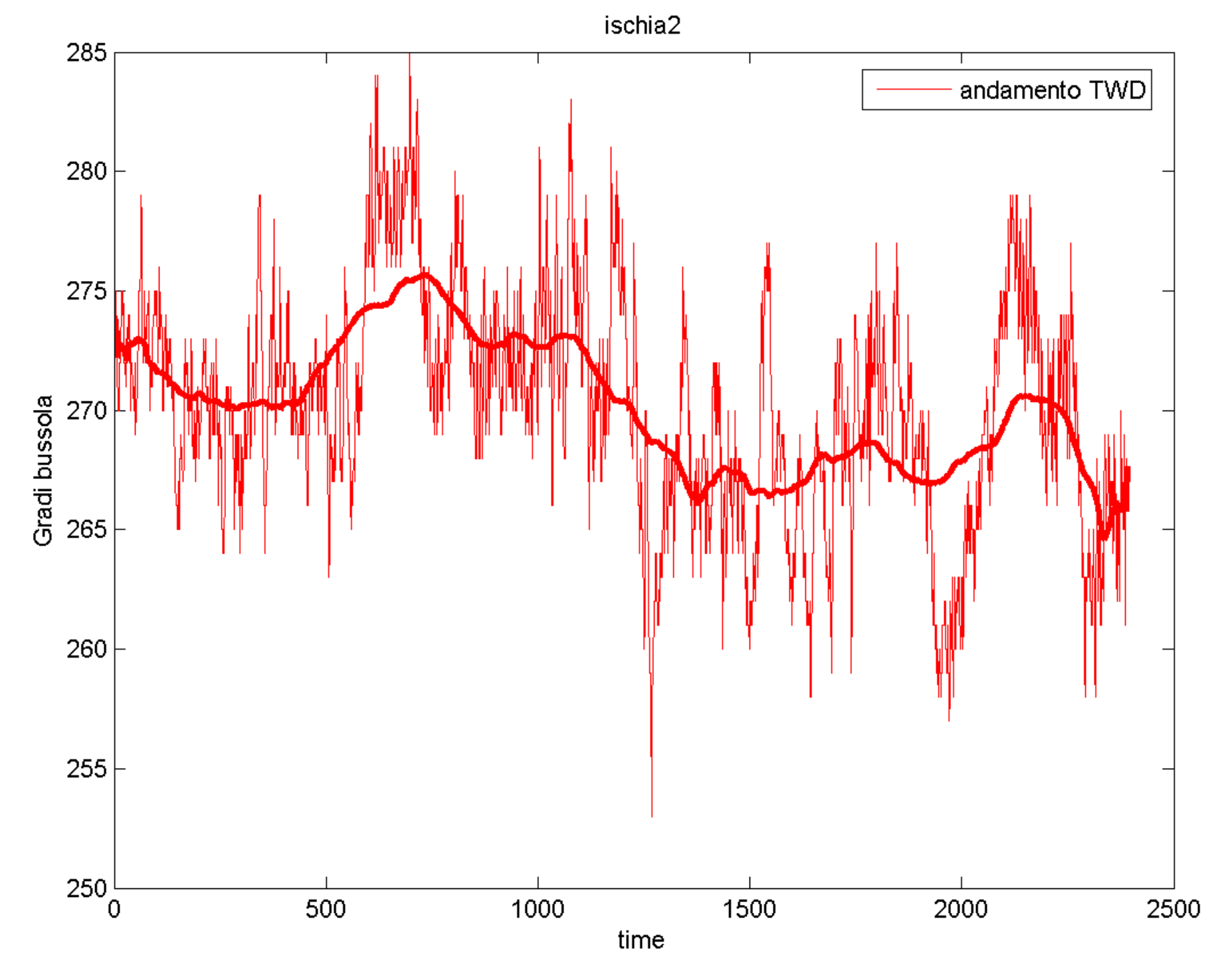}
\includegraphics[height=5.4cm]{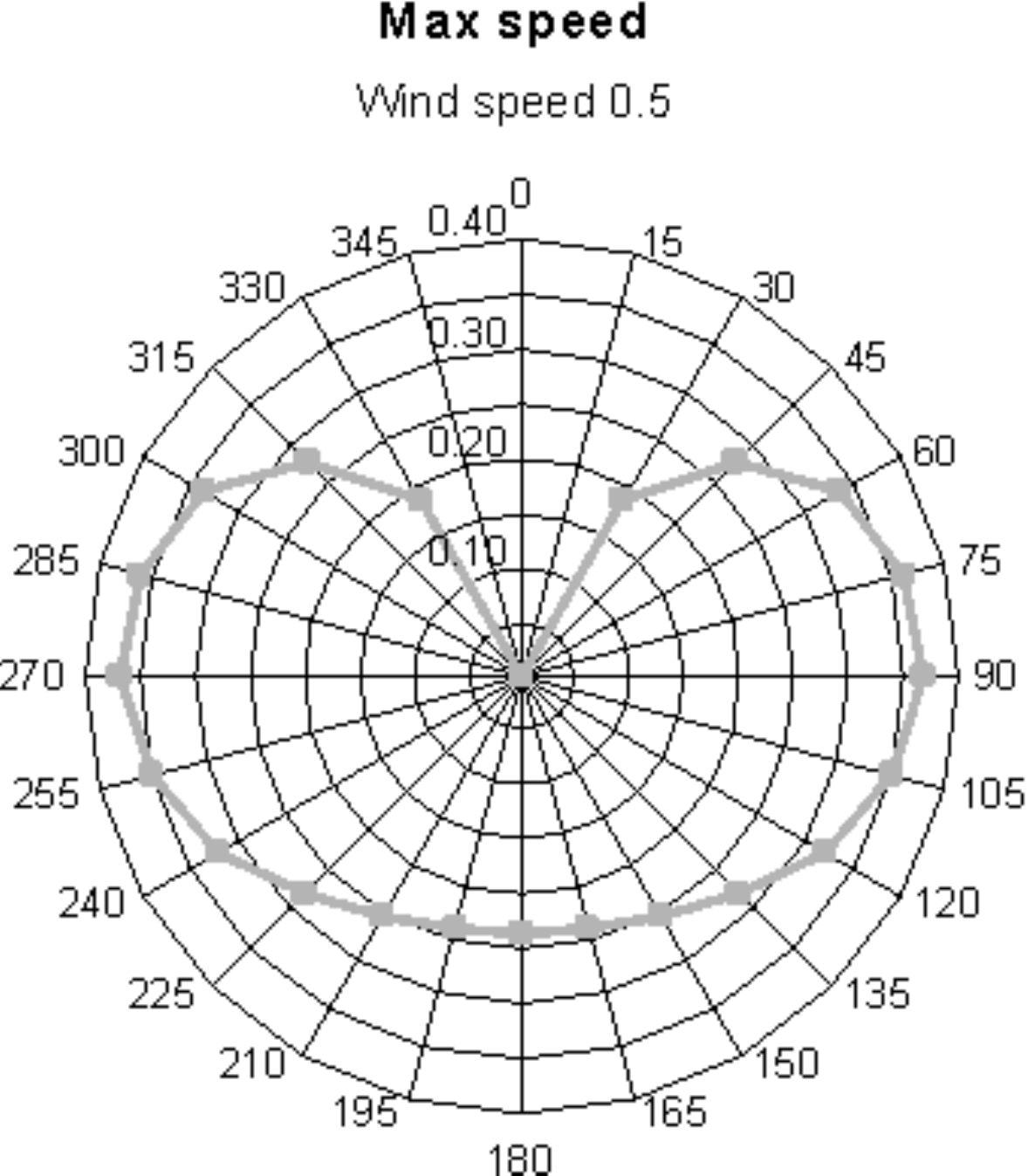}\hspace{1cm}
\caption{(left)  Evolution of the wind observed during a race (thickened plot the average over 180s) and Polar plot of the boat speed (in function of the TWA, with fixed TWS )(right)} \label{fig:1}
\end{center}
\end{figure}

In the most typical and basic form, the route planning problem requires reaching a windward mark in minimum time within a variable wind field. In particular, according to the previous literature and to experimental evidence, this problem is approached from a stochastic viewpoint, in which the wind field will be modeled as having both a deterministic and a stochastic component (cf. \cite{ZM}). In fact, even having a reliable forecast of its evolution, the wind field is affected by a random fluctuation (Fig. \ref{fig:1}, left) that represents a crucial part of the analysis.  Then, the optimal strategy requires using wind variations so as to minimize the arrival time.

Among the various techniques used to tackle the route planning problem, we will focus here on Dynamic Programming. A discrete Markov chain approximation of the stochastic route planning problem has been proposed, for example, in \cite{PM01}, based on an assumption of Markovian behavior for the wind model. We also quote the thesis \cite{T15}, which concentrates on the modeling and racing strategy for the specific case of match races, but also provides a nice review of the relevant literature.

The velocity of a sailing boat is usually characterized via the so-called \emph{polar plot} (see Fig. \ref{fig:1}, left), which shows, at fixed wind speed, the boat velocity as a function of the angle between the boat direction and the direction (True Wind Angle or TWA) and speed (True Wind Speed or TWS) of the wind. A common feature of polar plots is to present vanishing boat speeds when heading in the direction of the wind; this implies that, when a change of direction is performed passing through the zero-TWA (i.e., when {\em tacking}), the boat slows down. This delay is a key point in short-course racing: avoiding taking it into account may result in unrealistic optimal paths, possibly heading directly against the wind via an infinitely fast switching (``chattering'') between a positive and a negative angle.

A very natural form of taking into account the delay effects of tacking, without increasing the dimension of the problem, is to use the framework of {\em hybrid systems}. This notion has been proposed in the 90s to treat control systems operated by both continuous and discrete control actions. Among the extensive literature on hybrid control, we quote here the general approach proposed in \cite{branicky1998unified}, along with the study in \cite{MR2177312,BM97}, aimed at treating optimal hybrid control problems in the framework of viscosity solutions, for respectively the deterministic and the stochastic case (see also \cite{BCD97}). In the present work, this specific framework is used to obtain a sound theoretical background \cite{BS91} for the convergence of monotone schemes. Previous works in this direction are \cite{MR3328207,FS16}, both of which, however, treat the deterministic case. A stochastic route planning problem, close to the hybrid formulation although with a slightly different model, has been considered in \cite{SV15}.

In the problem under consideration, we will consider as continuous control the choice of a particular route {\em while keeping the tack unchanged}, and as discrete control a change of tack (and, possibly, any other discrete action like a change in the sail configuration). Moreover, we will also adopt a Dynamic Programming strategy, which, in the relatively low dimension of the problem, gives some definite advantages:

\begin{itemize}

\item we obtain a \emph{static feedback control} which is computed once and for all for a given target and wind model;

\item costs associated to discrete controls (tacking, changes of configuration) may be easily taken into account;

\item we can deal with the presence of a \emph{stochastic term} in the system, modeling the variations of the wind with respect to the expected data;

\item the presence of \emph{state constraints} (coasts, obstacles), as well as current and tides, can be handled in a relatively straightforward way.

\end{itemize}

The paper is structured as follows: in Section \ref{Sect:hyb} we introduce the mathematical framework of stochastic hybrid systems. In Section \ref{Sect:route} we discuss in detail the modeling of the route planning problem. Section \ref{Sect:scheme} describes and analyzes the numerical framework, and shows a monotone numerical scheme suitable for the problem under consideration. Some hints about the computational implementation are also provided. In Section \ref{Sect:test} we perform various tests to validate the approach and show the effectiveness of the technique in different scenarios of application.

\section{A hybrid framework for the route planning problem} \label{Sect:hyb}

Let us introduce the mathematical description of the problem. 
A {\em hybrid control system} is a system that, in addition to the usual continuous control, can undergo discrete control action, like switching between different dynamics or jumping from one state to another, in a discontinuous way. A sailing boat is an example of this kind of systems: with the aim of controlling the craft, the crew has the possibility to adjust the direction, switch between tacks (a boat keeping a TWA between $0^o$ and $180^o$ is told to be on \emph{port tack}, and on the \emph{starboard tack} otherwise), and change the number and the configuration of the sails. Any of these discrete choices will be associated to a different dynamics. Note that any change of configuration and/or tack implies a speed loss, which will come into play as soon as we will define a cost functional for this system.

Among the various mathematical formulations describing a hybrid system, we refer here to a simplified version of the one proposed in \cite{BM97} (similar formulations for the deterministic case have been proposed in \cite{branicky1998unified, MR2177312}, while a stochastic hybrid problem much in the same spirit of the one under consideration has been studied in \cite{F14}). Let $\I=\{1,2,\ldots,N_\I\}$ be finite, and consider the controlled system $(X,Q)$ described by:
\begin{equation} \label{eq_stato}
 \begin{cases}
 dX(t)=f(X(t),Q(t),u(t))dt+\sigma(X(t),Q(t))\,dW_t, \\
 X(0)= x, \ Q(0^+)=q,
 \end{cases}
\end{equation}
where $x,X\in \RR^d $, $q,Q\in \I$ and $dW_t$ is the differential of a $d$-dimensional standard Brownian process. Here, $X(t)$ and $Q(t)$ denote respectively the continuous and the discrete component of the state at time $t$, and, in order to end up with a stationary Dynamic Programming equation, we are assuming that $f$ depends on $t$ only via $X$, $Q$ and $u$. The function $f :\RR^d\times \I\times U \to \RR^d$ represents the continuous dynamics, for a set of continuous controls given by:
$$
\mathcal{U}=\{u:(0,\infty) \to U \> | \> u \text{ measurable}, \ U \mbox{ compact} \},
$$
and we assume both $f$ and $\sigma$ to be globally bounded and uniformly Lipschitz continuous w.r.t.~$x$.

The term $Q(t)$ models the possibility to switch between the various dynamics of the system, and takes values in the set of piecewise constant discrete controls $\cQ$, that is:
$$
\cQ = \{Q(\cdot):(0,\infty) \to \I \> | \> Q(t)=\sum_i^N w_i \chi_{t_i}(t) \},
$$
where $\chi_i(t)=1$ if $t\in [t_i,t_{i+1})$ and $0$ otherwise,  $\{t_i\}_{i=1,...,N}$ are the (ordered) times at which a switch occurs, and $\{w_i\}_{i=1,...,N}$ are values in $\I$. With respect to the more general setting, we are assuming some simplifications, and in particular:

\begin{itemize}

\item
The discrete control is in the form of a {\em switching}, i.e., it can only change the discrete component of the state $Q(t)$. In other terms, jumps in the $X$-component of the evolution are not allowed;

\item
A controlled switching can occur in the whole state set $\R^d\times\I\subseteq\RR^d\times\I$; on the other hand, we will avoid mandatory (``autonomous'') switchings. It is not difficult to consider a more general approach \cite{BM97}, in which the state takes values in a bounded set $\Omega$, and autonomous switchings (or a stopping cost) are imposed at the boundary of $\Omega$, if required (e.g., because of obstacles or state constraints);

\item
At this stage, nothing prevents the system from undergoing a ``Zeno effect'' (accumulation of switching times); however, the assumptions on the cost functional will preclude this situation.

\end{itemize}

The trajectory starts from $(x,q)\in \RR^d\times \I$. The choice of the control strategy defined as $\mathcal{S}:=\left(u,\{t_i\},\{Q(t_i^+)\}\right)$ has the objective of minimizing the following cost functional of minimum time type:
\begin{equation}\label{J}
 J(x,q;\mathcal{S})  := \E \left(\int_0^{\tau_{x,q}} e^{-\lambda t} dt + \sum_{i=0}^N C\left(X(t_i),Q(t_i^-),Q(t_i^+)\right)e^{-\lambda t_i} \right)
\end{equation}
where $\tau_{x,q}$ is the first time of arrival in a given compact target set $\cT\subset\RR^d$, i.e.,
$$
\tau_{x,q}:=\min_{t\in [0,+\infty)}\{t\;| \;X(t)\in\cT\},
$$
$\lambda>0$ is the discount factor, and $C:\RR^d\times \I \times \I  \to \RR_+$ is the switching cost between the dynamics, which is assumed to have a strictly positive infimum, to be bounded and Lipschitz continuous w.r.t.~$x$ and to satisfy the further condition
\begin{equation}\label{triang}
C(x,q_1,q_2) < C(x,q_1,q_3) + C(x,q_2,q_2),
\end{equation}
for any triple of indices $q_1$, $q_2$ and $q_3$.

The value function $v$ of the problem is then defined, for $\mathcal{S}\in\cU\times\R_+^{\N}\times\I^{\N}$, as:
\begin{eqnarray} \label{f_valore}
 v(x,q):= \inf_\mathcal{S} J(x,q;\mathcal{S}),
\end{eqnarray}
and is characterized via a suitable Hamilton--Jacobi--Bellman (HJB) equation. Continuity of the value function, which allows applying the framework of viscosity solutions, is a delicate matter in deterministic hybrid control problems (we refer the reader to \cite{MR2177312} for a precise set of assumptions), whereas in the stochastic case the literature reports somewhat weaker assumptions (see \cite{BM97}). Since the interest of this work is in a specific application, it might happen that some of the assumptions ensuring continuity will be dropped in what follows, but we will show nevertheless that the technique is effective and robust even in a more general setting.

We collect here a set of (not necessarily minimal) basic assumptions on the problem.

\subsubsection*{Basic assumptions}

\begin{enumerate}

\item
The functions $f$, $\sigma$ and $C$ are uniformly bounded, and Lipschitz continuous wrt $x$. The function $C$ has a positive infimum and satisfies \eqref{triang};

\item
The set of controlled switching coincides with the whole state space;

\item
$\lambda>0$.

\end{enumerate}

\bigskip

We recall some basic analytical results about the value function \eqref{f_valore}. Using a suitable generalization of the Dynamic Programming Principle it is possible to prove that the value function of the problem solves a Bellman equation in a Quasi-Variational Inequality form. More precisely, defining for $x,p\in\RR^d$ and $q\in\I$ the Hamiltonian function by
\begin{equation}\label{Ham}
 H(x,q,p) := \sup_{u\in U}\{ - f(x,q,u)\cdot p - 1 \}
\end{equation}
and the controlled switching operator $\cN$ by:
\begin{equation*}
 \cN \phi(x,q) := \inf_{w\in \I} \{\phi(x,w)+C(x,q,w)\},
\end{equation*}
we have a Bellman equation of the following form:
  \begin{equation}\label{hjb}
  \max\left(v-\cN v, \lambda v+ H(x,q,D v)+\frac{1}{2}\tr\left( \sigma\sigma^t D^2 v\right) \right)= 0,
  \end{equation}
defined on $(\R^d\setminus\cT)\times\I$, i.e., a system of Quasi-Variational Inequalities, complemented with the boundary condition
\[
v(x,q)=0 \quad (x\in\partial\cT).
\]
In what follows, we will assume that the problem is posed on a set $\Omega$ of the state space, and either that $\Omega=\RR^d$, or that a stopping cost is defined on the boundary of $\Omega$ (in the form of a weak Dirichlet condition, see \cite{falconeferretilibro}), so as to enforce state constraints by penalization, and work on a finite computational domain. We will give suitable examples in the numerical test section.

In \eqref{hjb}, we can identify two separate Bellman operators, which provide respectively the best possible switching, and the best possible continuous control. The argument attaining the maximum in \eqref{hjb} represents the overall optimal control strategy.

\section{Practical models for route planning}\label{Sect:route}

In order to apply the techniques introduced above to solve some problems of route planning, we outline in this section some general ideas towards a formal modeling of the problem. The idea is to start with some easy cases, useful to show the validity of the solutions produced, and to pass to more complex situations, proving that such approach can deal with many of the delicate points of the problem. 

In general, we expect that a reliable modeling should use at least a 3-dimensional state space, i.e., $d\ge 3$, except for the case in which the average direction of the wind is constant and its variations come only from the stochastic part (this will be cleared up in the forthcoming examples). In this setting, two components $x_1$ and $x_2$ of the state space represent the position of the boat, while the third component $x_3$ accounts for the evolution of the wind. Note that, in a different formulation, the third variable could be given by time, and in this case the wind should be provided as a vector field defined on $\RR^2\times [0,T]$, and the HJB equation would be time-dependent. However, we will rather pursue here a model leading to the stationary HJB equation \eqref{hjb}.

\begin{figure}
\centering
\includegraphics[height=6cm]{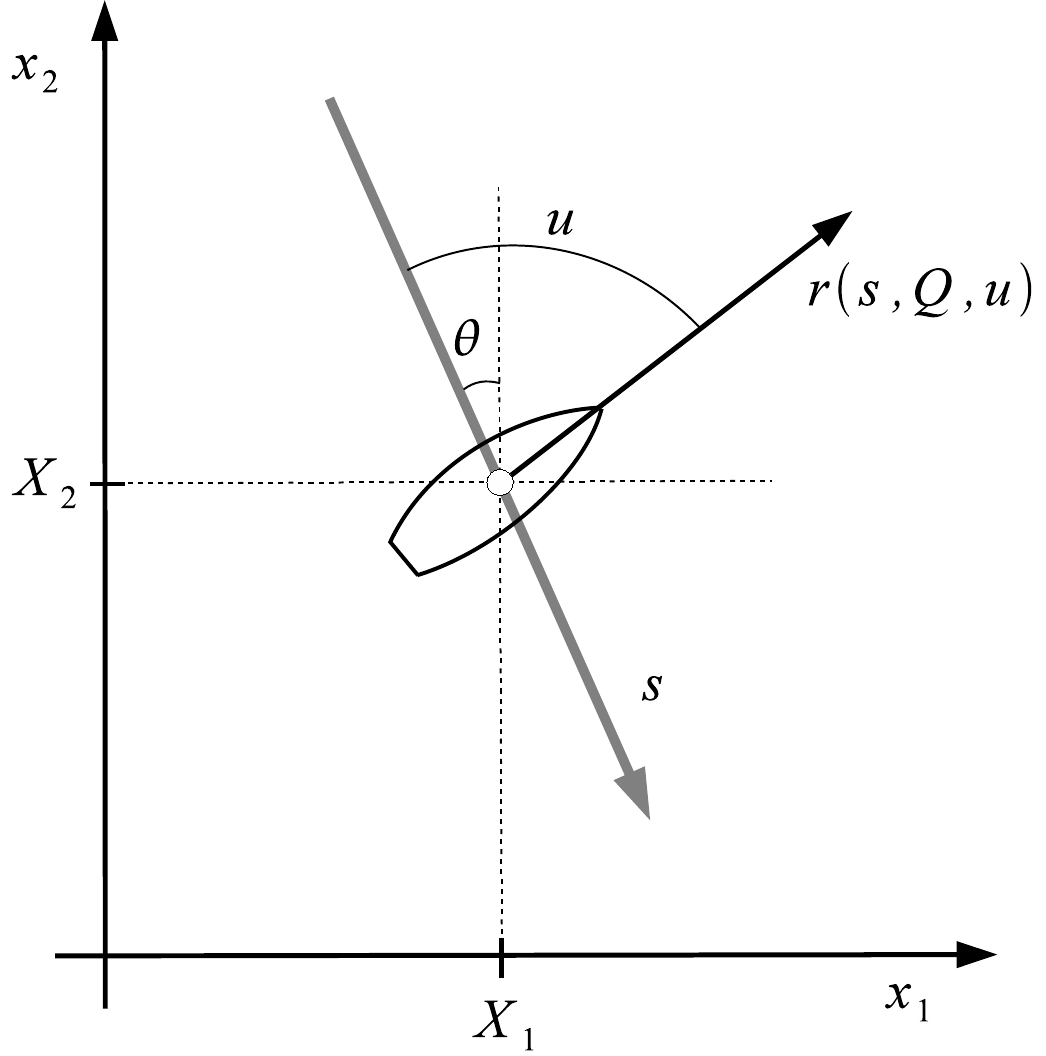}
\hspace{0.8cm}
\includegraphics[height=6cm]{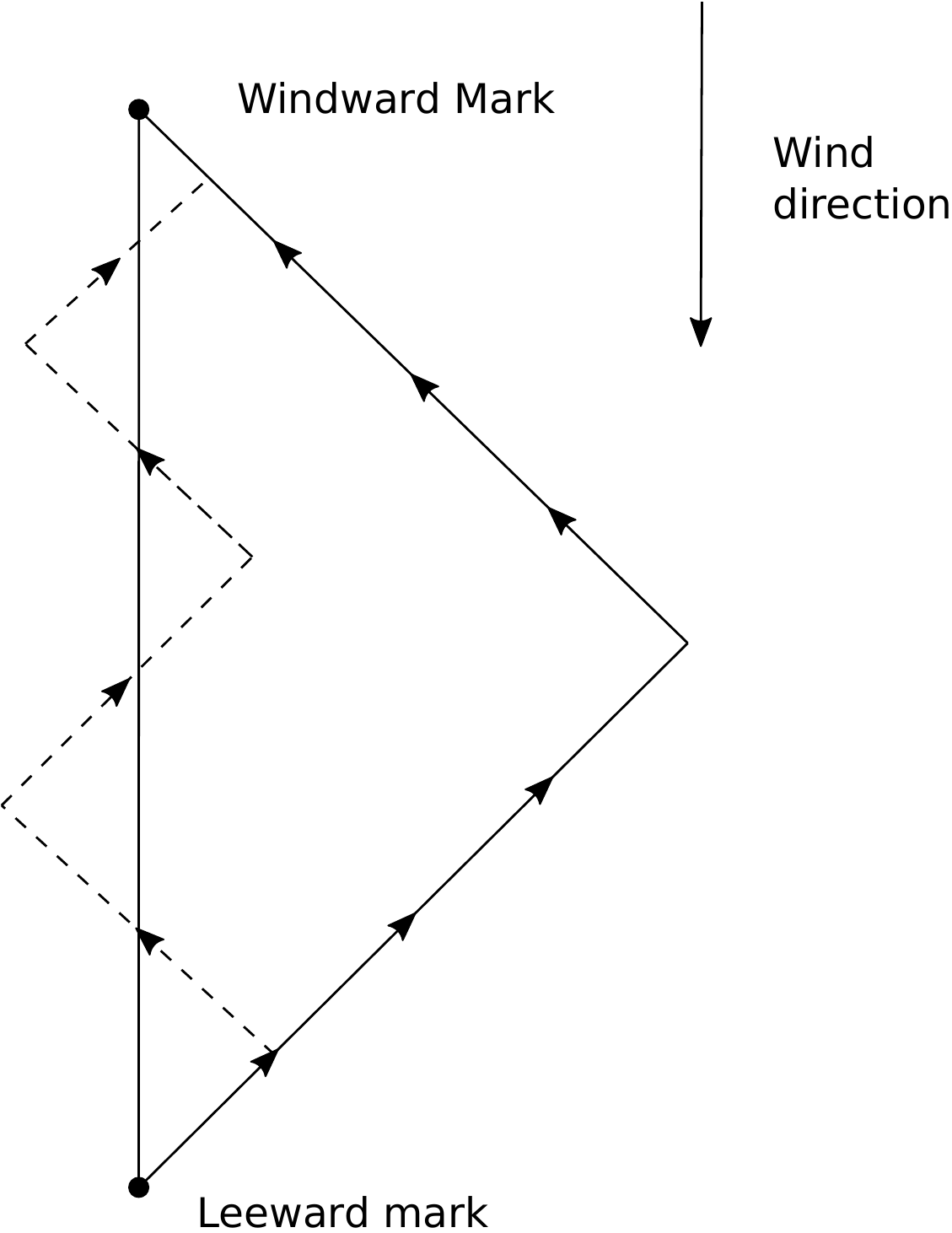}
\caption{(Left) Geometric setting: Wind (grey arrow), boat speed (black arrow) and control $u$, (right) upwind tacking: to go from the Leeward to the Windward mark it is not possible to point directly to the target.} \label{fig:notazioni}
\end{figure}

We start by describing the motion of the boat as resulting from both the wind vector field and the boat characteristics. Although more general forms could be used, we will make the standing assumption that the control $u$ denotes the (unsigned) angle between the boat direction and the wind, so that $u(t) \in U=[0,\pi]$, and that the stochastic component of the dynamics appears only in the wind evolution. The speed $r$ of the boat will be assumed to depend only on the wind speed $s$, the wind direction $\theta$ and the angle $u$, along with possible discrete controls described by $Q$. The motion of the boat is then described by
\begin{equation}\label{eq:motion}
\begin{cases}
\dot X_1(t) = r(s(X(t),t),Q(t),u(t))\sin(-\theta(X(t),t)\pm u(t)) \\
\dot X_2(t) = r(s(X(t),t),Q(t),u(t))\cos(\theta(X(t),t)\pm u(t)),
 \end{cases}
\end{equation}
where the plus sign corresponds to the starboard tack, and the minus to the port tack. The function $r:\R_+\times\I\times [0,\pi] \rightarrow \R_+$ models the \emph{polar plot} of the boat, and provides the boat speed as a function of the angle $u$ of the trajectory w.r.t.~the wind, of the wind speed $s$ and of the sail configuration. Note that with control $u=0$ the trajectory points directly in the `upwind' direction whereas if $u=\pi$ the trajectory has the same direction of the wind field. Fig. \ref{fig:notazioni} summarizes the geometric setting.

The choice of the function $r(s,q,u)$ is related to the technical characteristics of the craft. It differs from one boat to another, but we can detect some general features:

\begin{itemize}

\item
$r(s,q,u)$ is continuous w.r.t.~both variables $s$ and $u$;

\item
For given $q$, it depends on the wind speed $s$ and the relative angle $u$, but neither on time nor on position; moreover, the dependence on $s$ is monotone;

\item
$r(s,q,0)=0$, which means that the boat has always zero speed when `pointing directly against the wind';

\item
$r(s,q,\cdot)$ has typically (but not necessarily) a single maximum point inside $[0,\pi]$.

\end{itemize}

In practice, the function $r$ is determined via experimental measures. We will not try to make it more explicit here, and rather assume that it is defined by some polynomial interpolation of the experimental data.

The wind will be characterized by two functions, direction $\theta$ and speed $s$. These function may depend on position and time, but we assume that they evolve in time according to a lumped parameter model, i.e., the system of stochastic differential equations (SDEs)
\begin{equation}\label{eq:wind}
\left\{
\begin{array}{l}
ds(x,t) = g_1(x,s(x,t),\theta(x,t))dt + g_2(x,s(x,t),\theta(x,t))dW_t^{(1)} \\
d\theta(x,t) = h_1(x,s(x,t),\theta(x,t))dt + h_2(x,s(x,t),\theta(x,t))dW_t^{(2)},
\end{array}\right.
\end{equation}
in which $x$ is considered as a parameter, $dW_t^{(i)}$ $(i=1,2)$ denotes the differential of a standard Brownian process, and, depending on the complexity of the model, both $s$ and $\theta$ could be multi-dimensional. 

\subsection{A detailed example}\label{exampledetailed}

We will discuss now in detail a simplified model, and later examine possible generalizations. We assume that $d=3$, the third state variable being the wind direction, $x_3=\theta$, and, for simplicity, that $\I=\{1,2\}$, thus meaning that the discrete control consists only in tacking, and the function $r$ does not depend on $Q$. If, for example, the port tack is associated to the discrete state $q=1$ and the starboard tack to $q=2$, then we can write
\begin{equation}\label{eq:simpl_dyn}
\begin{cases}
\dot X_1(t) = r(s,u)\sin\left(-\theta + (-1)^{Q(t)} u\right) \\
\dot X_2(t) = r(s,u)\cos\left(\theta + (-1)^{Q(t)} u\right).
 \end{cases}
\end{equation}

\subsubsection*{Wind evolution}

A simple model, which still catches most features of route planning problems, at least in an upwind part of a race, is that in which the wind speed $s$ is kept constant w.r.t. position and time,
\[
s(x,t) \equiv \bar s,
\]
the direction is constant w.r.t.~position, and evolves according to the one-dimensional SDE
\begin{equation}\label{eq:simpl_wind}
d\theta = a(\theta)dt + \bar\sigma dW_t.
\end{equation}
In particular, taking $a(\theta)=\bar a$ (constant) would model a situation in which the wind direction comes from the superposition of a constant (clockwise or anti-clockwise) drift with a random fluctuation. In this case,
\[
\sigma = \begin{pmatrix}
0 \\
0 \\
\bar\sigma
\end{pmatrix},
\]
and the diffusion term in the HJB equation \eqref{hjb} is degenerate and acts only along the third dimension. Note that, roughly speaking, assuming a one-dimensional model for $\theta$ parallels the assumption of discrete Markov process without memory as in \cite{PM01}.

\subsubsection*{Simplified dynamics for a windward leg}

Since in a windward leg the boat is usually kept at the most efficient angle $u^*$ w.r.t.~the wind (typically, $u^*\approx \pi/4$), we can give a simpler dynamics for this case by assuming that the control $u$ is frozen at $u\equiv u^*$ and the only control action consists in tacking. This model might suffer from a lack of controllability with respect to the complete model; nevertheless, and despite its simplicity, it still provides qualitatively good results.

\begin{figure}
\centering
\includegraphics[height=4.7cm]{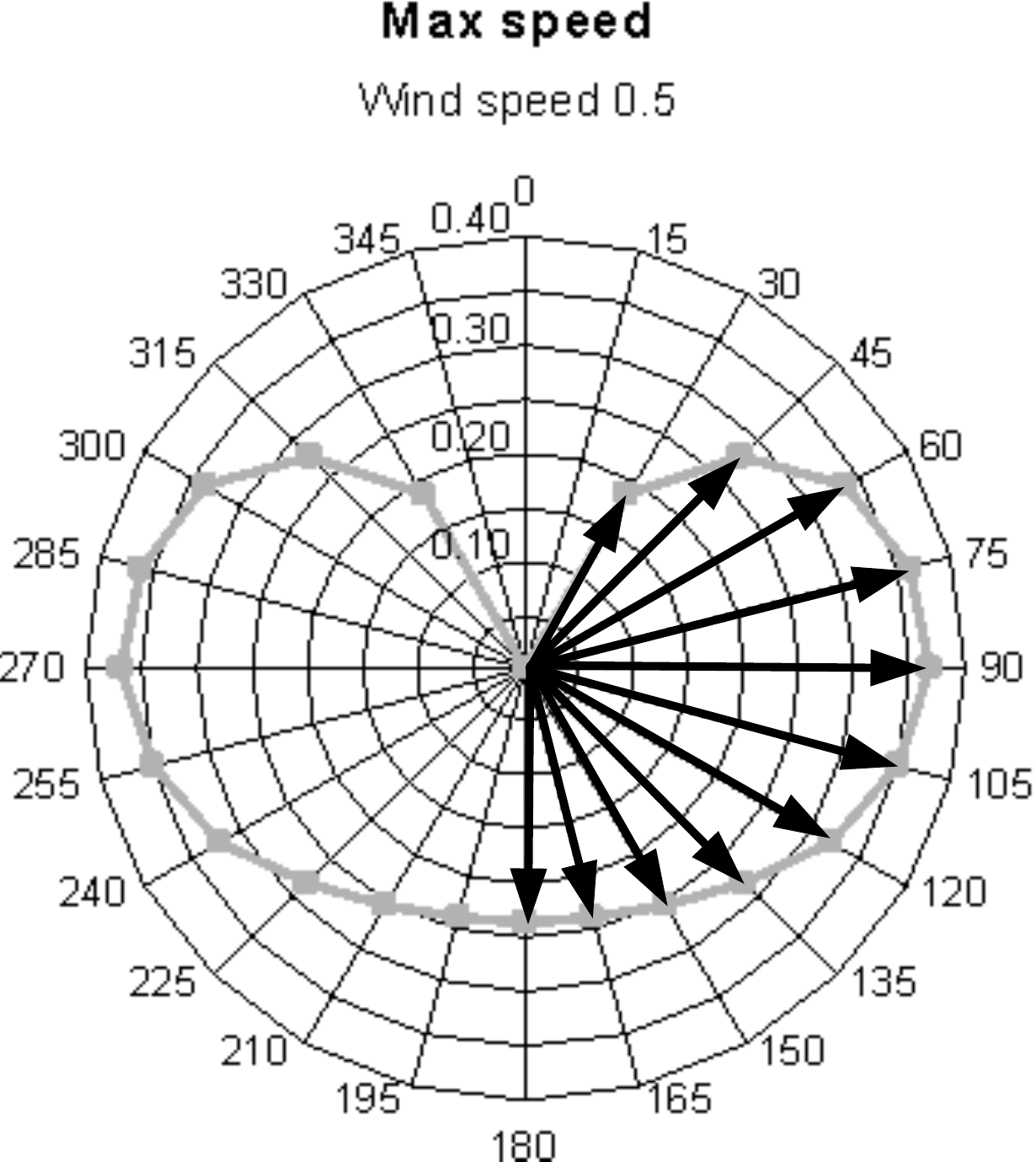}\hspace{.5cm}
\includegraphics[height=4.7cm]{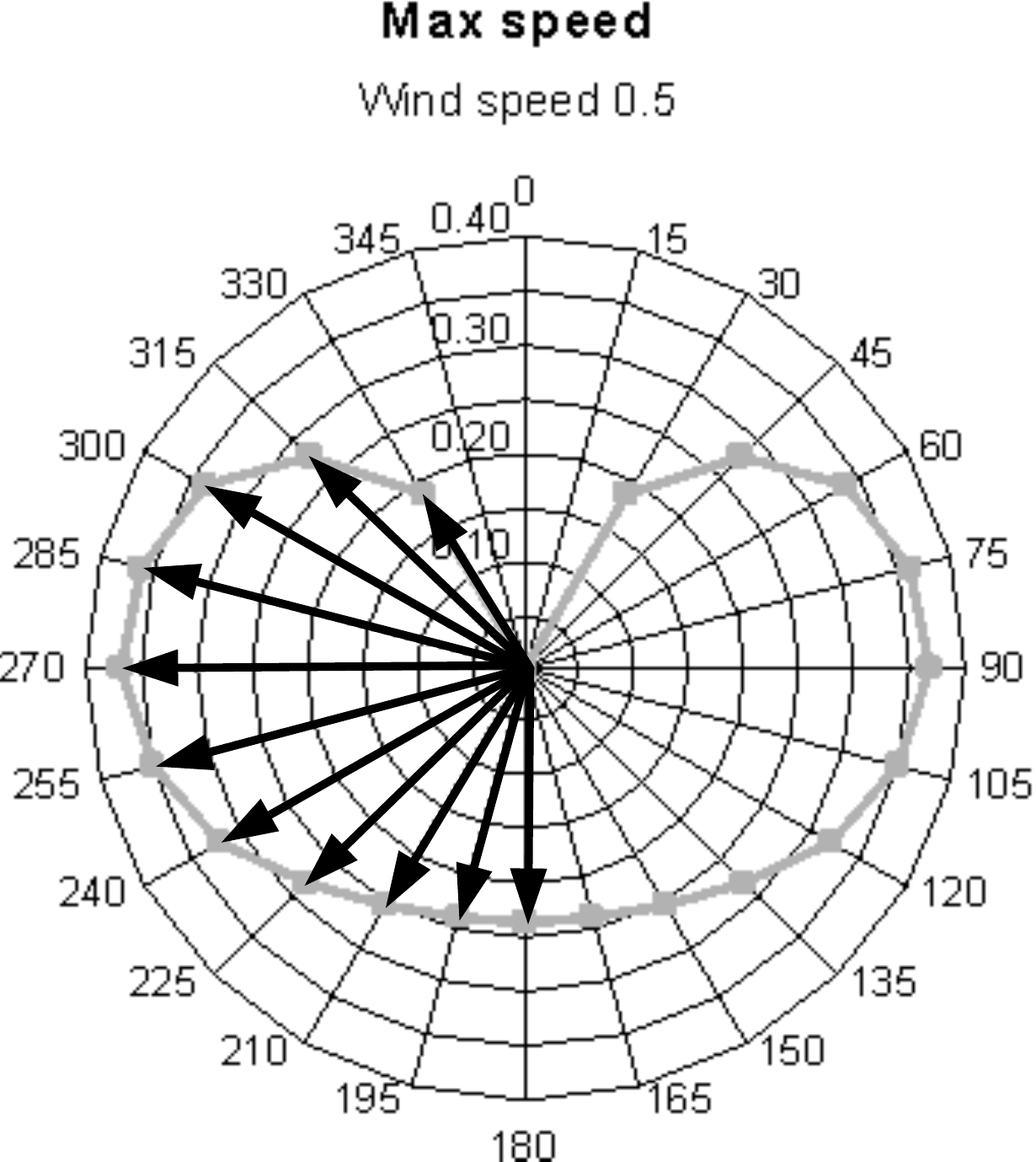} \\
complete \\
\vspace{1cm}
\includegraphics[height=4.7cm]{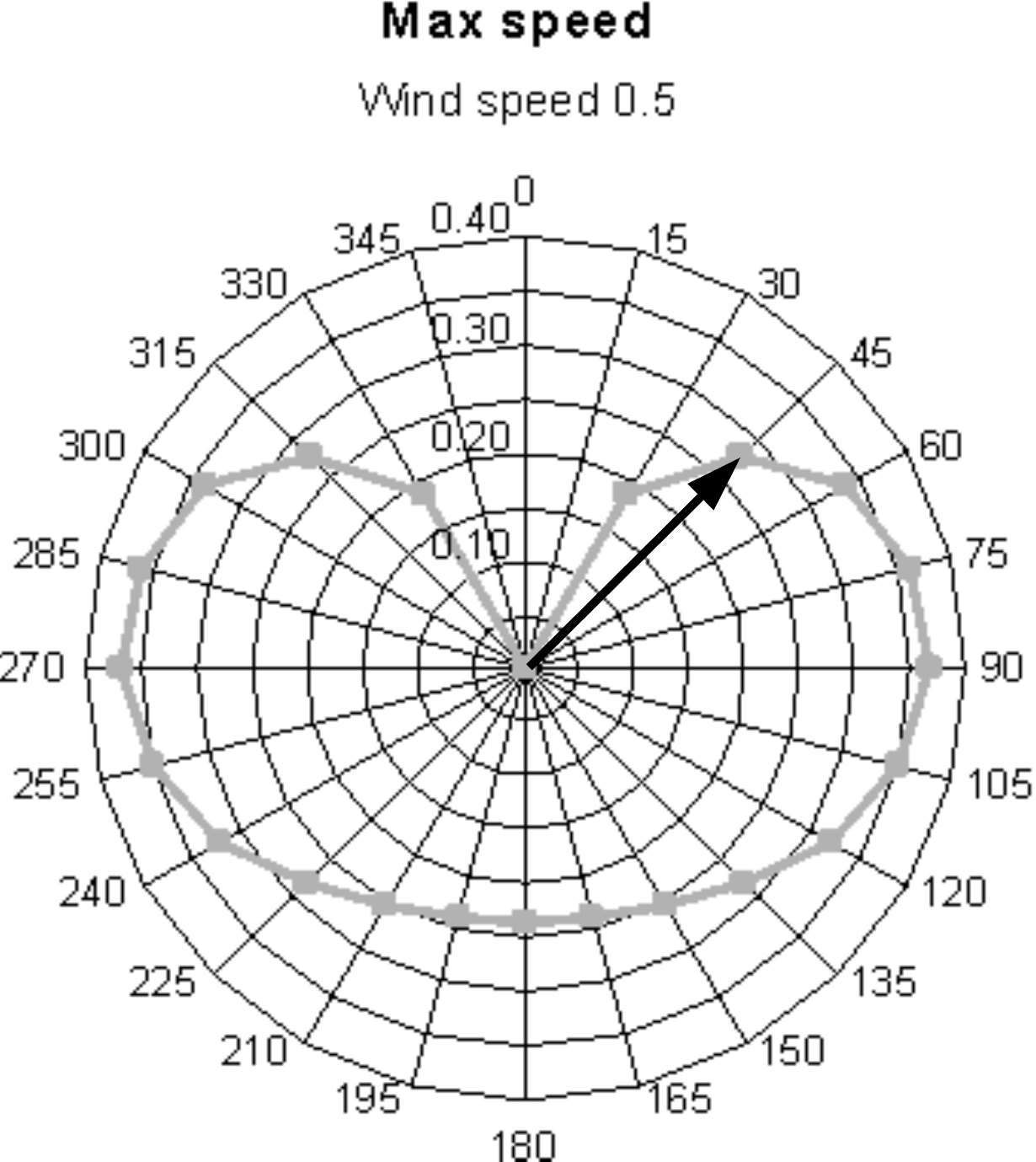}\hspace{.5cm}
\includegraphics[height=4.7cm]{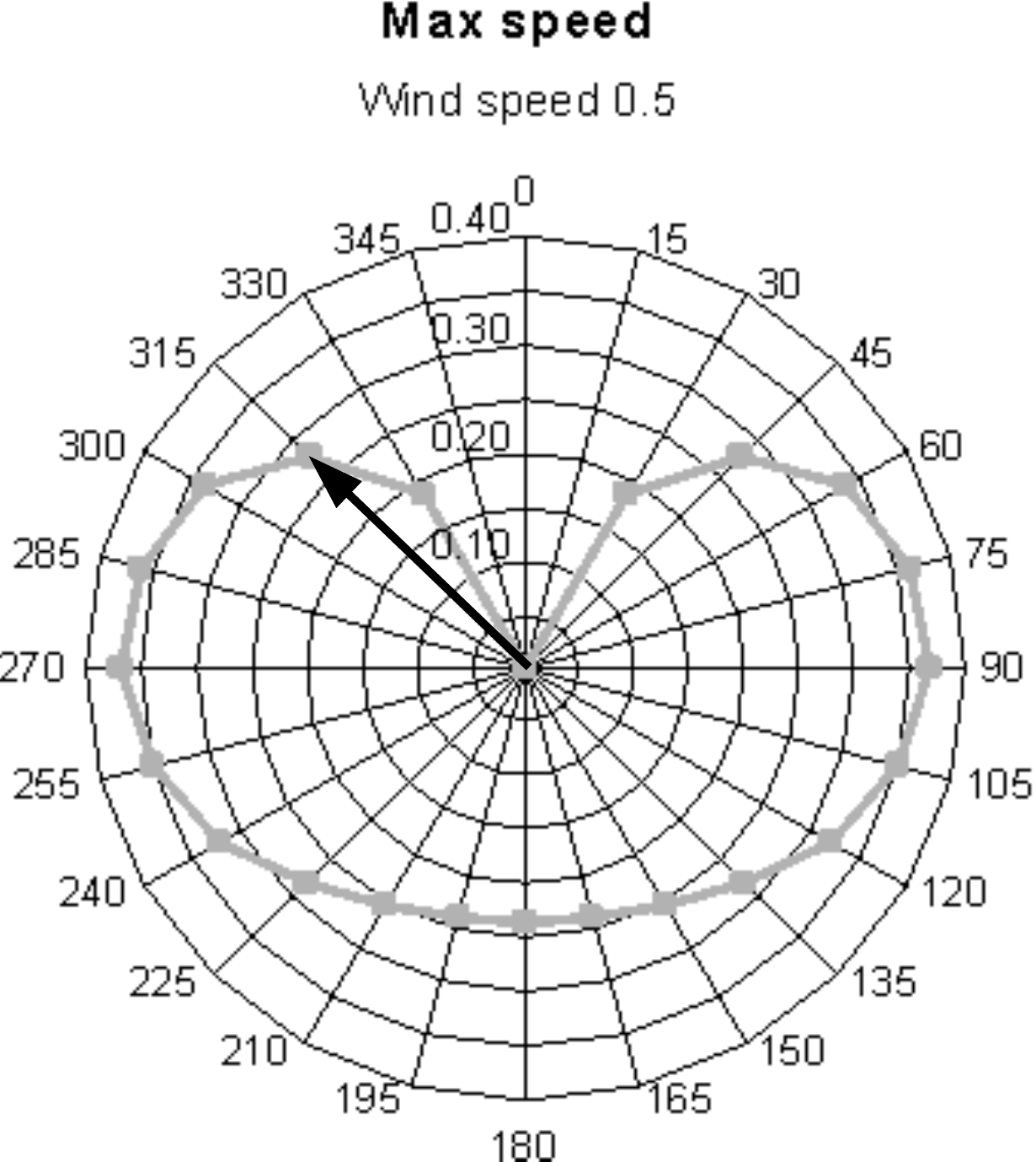} \\
simplified \\
$Q=1$ \hspace{3.5cm} $Q=2$
\caption{(upper) The two dynamics of the system, superposed on the polar plot of the speed; (lower) simplified dynamics based on the angle of largest windward component of the speed.} \label{fig:dyn}
\end{figure}

Fig. \ref{fig:dyn} compares complete and simplified dynamics.

\subsubsection*{Cost functional}

The cost functional is already in the specific form needed. We only remark that we can reasonably assume that the cost of tacking is constant, $C(x,1,2)=C(x,2,1)=\bar C$, and use a positive $\lambda$ in order to obtain a contractive Bellman operator. Note that, although strictly speaking a minimum time problem should not be discounted, yet the introduction of a discount factor may be seen as the application of the so-called {\em Kru\v{z}kov transform} (see \cite{falconeferretilibro}).

\subsubsection*{Resulting HJB equation}

We collect now all the information in an explicit HJB equation. Concerning the deterministic component of the dynamics we have, using \eqref{eq:simpl_wind} into \eqref{eq:simpl_dyn}:
\begin{equation*}
\begin{cases}
f_1(x,q,u) = r(\bar s,u)\sin\left(-x_3 + (-1)^q u\right) \\
f_2(x,q,u) = r(\bar s,u)\cos\left(x_3 + (-1)^q u\right) \\
f_3(x,q,u) = a(x_3),
 \end{cases}
\end{equation*}
which results in a Hamiltonian function of the form
\begin{multline*}
 H(x,q,p) = \sup_{u\in U}\big\{ - r(\bar s,u)\left[\sin\left(-x_3 +(-1)^q u\right)p_1 + \cos\left(x_3 +(-1)^q u\right)p_2\right] \big\}\\ - a(x_3)p_3 - 1
\end{multline*}
if a continuous control action $u$ is possible, and in the form
\begin{equation*}
 H(x,q,p) = - r(\bar s,u^*)\left[\sin\left(-x_3 +(-1)^q u^*\right)p_1 + \cos\left(x_3 +(-1)^q u^*\right)p_2\right] - a(x_3)p_3 - 1
\end{equation*}
if $u\equiv u^*$ and the only control action is to tack (note that, in this case, the Hamiltonian function is linear in $p$). The term related to the stochastic component of the evolution reads in turn as
$$
\frac{1}{2}\tr\left( \sigma\sigma^t D^2 v\right) = \frac{\bar\sigma^2}{2}v_{x_3x_3}.
$$
Last, assuming a constant switching cost $\bar C$, the switching operator $\cN$ is given by:
\begin{equation*}
 \cN \phi(x,q) = 
 \begin{cases}
 \min \{\phi(x,1), \phi(x,2)+\bar C\} & \text{ if } q=1, \\
 \min \{\phi(x,2), \phi(x,1)+\bar C\} & \text{ if } q=2.
 \end{cases}
\end{equation*}
Then, the Bellman equation takes the form of the following system of two Quasi-Variational Inequalities:
  \begin{multline} \label{qvi1}
 \max\Big(v(x,1) -\min \left\{v(x,1), v(x,2)+\bar C\right\}\>,\> \lambda v(x,1) \\
 +\sup_{u\in U} \big\{- r(\bar s,u) \left[\sin\left(-x_3 +(-1)^q u\right)v_{x_1}(x,1) + \cos\left(x_3 +(-1)^q u\right)v_{x_2}(x,1)\right]\big\} \\
- a(x_3)v_{x_3}(x,1) - 1+\frac{\bar\sigma}{2}v_{x_3x_3}(x,1) \Big)= 0, 
\end{multline}
  \begin{multline} \label{qvi2}
  \max\Big(v(x,2) -\min \left\{v(x,2), v(x,1)+\bar C\right\}\>,\> \lambda v(x,2) \\
 +\sup_{u\in U} \big\{- r(\bar s,u) \left[\sin\left(-x_3 +(-1)^q u\right)v_{x_1}(x,2) + \cos\left(x_3 +(-1)^q u\right)v_{x_2}(x,2)\right]\big\} \\
 - a(x_3)v_{x_3}(x,2) - 1+\frac{\bar\sigma}{2}v_{x_3x_3}(x,2) \Big)= 0,
\end{multline}
in which the $\sup_u$ is dropped in the simplified case $u\equiv u^*$.

\subsection{Some remarks on a possible hierarchy of models}

While the previous model seems a reasonably simple way of taking into account the major processes involved in the route planning problem, we sketch some related models which lead to a different definition of the state space.

\begin{itemize}

\item
{\it Purely brownian wind direction/speed.} This is the only case leading to a lower dimension of the problem -- in this case, $d=2$. Assume that the wind is given by \eqref{eq:simpl_wind}, with $a(\theta)\equiv 0$. By a linearization of \eqref{eq:simpl_dyn} with respect to $\theta$ around the direction $\theta=0$, we can write an approximate dynamics in the form
\begin{equation*}
\begin{cases}
d X_1(t) = r(s,u)\sin\left((-1)^{Q(t)} u\right) dt - r(s,u)\cos\left((-1)^{Q(t)} u\right)\bar\sigma dW_t \\
d X_2(t) = r(s,u)\cos\left((-1)^{Q(t)} u\right) dt - r(s,u)\sin\left((-1)^{Q(t)} u\right)\bar\sigma dW_t,
 \end{cases}
\end{equation*}
in which the dependence upon the state variable $x_3$ has disappeared. A brownian evolution of the wind speed might be treated accordingly, by a linearization of $r(s,u)$ with respect to $s$.

\item {\it $x$-dependence of the wind direction/speed.} Note that it is possible to introduce in \eqref{eq:simpl_dyn} and \eqref{eq:simpl_wind} an $x$-dependence without increasing the conceptual difficulty of the problem. This change only results in making the Hamiltonian depend on $x_1$ and $x_2$ and has no consequences on the definition of the state space.

\item {\it 1D (possibly stochastic) evolution for wind speed.}In this case, the dynamics \eqref{eq:simpl_dyn} for the wind would be replaced by
\begin{equation*}
\begin{cases}
d\theta = a(\theta,s)dt + \bar\sigma_1 dW_t^{(1)} \\
ds = b(\theta,s)dt + \bar\sigma_2 dW_t^{(2)}.
\end{cases}
\end{equation*}
As a result, we would define $x_4=s$ and obtain a state space of dimension $d=4$.

\item
{\it More complex models for $\theta$ and/or $s$.} By a straightforward extension of the previous arguments, it is clear that the use of a $K_\theta$-dimensional model for $\theta$ and a $K_s$-dimensional model for $s$ results in a state space of dimension $d=2+K_\theta+K_s$.

\item
{\it Different configurations of the sails.} If $N_c$ different configurations of the sails are possible, then a change of configuration should appear as a discrete control. Since every sail configuration should be replicated on both the starboard and the port tack, we would then obtain $N_\I=2N_c$, with unchanged dimension $d$.

\item
{\it Long courses racing.} When long course racing is concerned, a local wind model like \eqref{eq:simpl_wind} is no longer possible. In this case, the deterministic component in the wind model might be dependent on time (and could be taken, e.g., from suitable weather forecasts like GRIB data), and in the simplest case the problem could be modelled via two state variables plus time. Accordingly, the HJB equation would be recast in the time-dependent form, with $d=2$, and $N_\I=2N_c$.

\end{itemize}

\section{Numerical solution via monotone schemes} \label{Sect:scheme}

In order to set up a numerical approximation for \eqref{hjb}, we construct a discrete grid of nodes $(x_j,q)$ in the state space with discretization parameters $\Delta x$ and, possibly, $\Delta t$ (the time discretization should be understood in the ``time marching'' sense). In what follows, we denote the discretization steps in compact form by $\Delta=(\Delta t, \Delta x)$ and the approximate value function by $V^\Delta$.

Following \cite{MR3328207}, we write the scheme at $(x_j,q)$ in fixed point form as
\begin{equation}\label{scheme1}
V^\Delta(x_j,q) = \min \left( NV^\Delta(x_j,q), \Sigma\left(x_j,q,V^\Delta\right)\right).
\end{equation}
In \eqref{scheme1}, the numerical operator $\Sigma$ is related to the continuous control, or, in other terms, to the approximation of the Hamiltonian function \eqref{Ham}. The discrete switch operator $N$ is computed at a node $(x_j,q)$ as
\begin{equation}\label{scheme11}
NV^\Delta(x_j,q) := \min_{q'\in \I}\left\{V^\Delta(x_j,q') + C(x_j,q,q')\right\},
\end{equation}
and, in fact, this corresponds to the exact definition.

\subsection{Theoretical analysis for monotone schemes}

The theoretical analysis of monotone schemes for HJB equations arising in sto\-cha\-stic hybrid control problems can be carried out using (with some adaptations) the arguments in \cite{MR3328207}. We start by proving that the value iteration for \eqref{scheme1}, i.e.,
\begin{equation}\label{scheme1_vi}
V^\Delta_{k+1}(x_j,q) = \min \left( NV^\Delta_k(x_j,q), \Sigma\left(x_j,q,V^\Delta_k\right)\right)
\end{equation}
is convergent. To this end, we denote by $S(V^\Delta_k)$ the vector (indexed by $j$) of the values at the right-hand side of \eqref{scheme1_vi}.

\begin{enumerate}

\item {\em If the scheme $\Sigma$ is monotone, then the mapping $S(\cdot)$ is also monotone.}

In fact, assuming that $U\ge W$, the inequality should be understood element by element. Since, as proved in \cite{MR3328207}, the operator $N$ is monotone, we only have to prove that the min of two monotone operators is monotone, i.e., that $S(U)\ge S(W)$. We have:
\begin{equation*}\label{mon1}
S(U) - S(W) = \min \left( NU(x_j,q), \Sigma\left(x_j,q,U\right)\right) - \min \left( NW(x_j,q), \Sigma\left(x_j,q,W\right)\right).
\end{equation*}
If in the two terms at the right-hand side of the equation above the minimum is achieved by the same operator, we have nothing else to prove, each of the two being monotone. Otherwise, if for example $\min \left( NU(x_j,q), \Sigma\left(x_j,q,U\right)\right)=NU(x_j,q)$ and $\min \left( NW(x_j,q), \Sigma\left(x_j,q,W\right)\right)=\Sigma\left(x_j,q,W\right)$ then
\begin{eqnarray*}
S(U) - S(W) & = & NU(x_j,q) - \Sigma\left(x_j,q,W\right) \\
& \ge & NU(x_j,q) - NW(x_j,q) \ge 0 
\end{eqnarray*}
and $S$ satisfies again the monotonicity condition. The same idea applies to the opposite case.

\item {\em The sequence $V^\Delta_k$ is monotone decreasing if $V^\Delta_0$ is chosen so that $V^\Delta_0 \ge S(V^\Delta_0)$.}

In fact, this condition implies that $V^\Delta_1 \le V^\Delta_0$, so that by monotonicity of $S$ we have $S(V^\Delta_1) \le S(V^\Delta_0)$, that is, $V^\Delta_2 \le V^\Delta_1$, and, inductively, $V^\Delta_{k+1} \le V^\Delta_k$. Note that this also entails that the scheme is stable in the $\infty$-norm.

\item {\em If, for any admissible $\Delta$, $V^\Delta_k$ is positive, then $V^\Delta_k$ converges to a fixed point $V^\Delta$ solution of \eqref{scheme1}.}

We recall that, under the assumptions made on the cost functional, positivity of solutions is a natural assumption, which is typically satisfied by monotone schemes (Upwind, Lax--Friedrichs, Semi-Lagrangian) under stability conditions. Convergence to a fixed point follows then from monotonicity and boundedness for all elements of the vectors $V^\Delta_k$.

\end{enumerate}

Last, once ensured that \eqref{scheme1} has a solution (which can be obtained by value iteration), the convergence of $V^\Delta$ to the value function $v$ is ensured by the Barles--Souganidis theorem \cite{BS91}, provided the scheme $\Sigma$ is consistent with the Hamiltonian function \eqref{Ham}, since a comparison principle holds for the continuous problem (see \cite{BM97}). We have therefore the following

\begin{theorem}
Let the basic assumptions hold. Let moreover $V^\Delta_0$ be chosen so that $V^\Delta_0 \ge S(V^\Delta_0)$. If the operator $\Sigma$ is monotone and such that $S(V)\ge 0$ for any $V\ge 0$, then the value iteration \eqref{scheme1_vi} converges to $V^\Delta$ solution of \eqref{scheme1}. Moreover, if the scheme $V=\Sigma(V)$ is also consistent with the equation
\[
\lambda v+ H(x,q,D v)+\frac{1}{2}\tr\left( \sigma\sigma^t D^2 v\right) = 0,
\]
then $V^\Delta(x_j,q)\to v(x_j,q)$ for $\Delta\to 0$.
\end{theorem}

\subsection{Example: a Semi-Lagrangian scheme}

A viable technique for solving \eqref{hjb} is a semi-Lagrangian scheme, obtained by adapting the scheme proposed in \cite{CF95}. The main advantage of such approach is an unconditional stability of the scheme with respect to the discretization parameters, still keeping monotonicity.

The scheme requires extending the node values to all $x\in\RR^d$ using an interpolation $\mathbb{I}$. We denote by $\mathbb{I}\left[V^\Delta\right] (x,q)$ the interpolation of the values $V^\Delta(x_j,q)$ computed at $(x,q)$. With this notation, a standard semi-Lagrangian discretization of the Hamiltonian and the diffusive term is given (see \cite{CF95}) by
\begin{equation}\label{scheme2}
\Sigma\left(x_j,q,V^\Delta\right)= \Delta t + e^{-\lambda \Delta t}\min_{u\in U}\left\{\frac{1}{2d}\sum_{i=1}^d\left(\mathbb{I}\left[V^\Delta\right] (\delta^+_i, q)+\mathbb{I}\left[V^\Delta\right] (\delta^-_i, q)\right)\right\},
\end{equation}
where
\[
\delta^\pm_i:= x_j+\Delta t \> f(x_j,q,u)\pm\sqrt{d \Delta t }\> \sigma(x_j,q)e_i.
\]
The full scheme is obtained by using \eqref{scheme11}--\eqref{scheme2} in \eqref{scheme1}.

As far as the interpolation $\mathbb I$ is monotone, the resulting scheme is consistent, monotone and $L^\infty$ stable, and therefore convergent via Barles--Souganidis theorem. Typical  monotone examples are $\PP_1$ (piecewise linear on triangles/tetrahedra) and $\Q_1$ (piecewise multilinear on rectangles) interpolations.

\subsection{Acceleration techniques}

Note that, in the basic setting, the use of a positive discount factor $\lambda$ ensures convergence of the value iteration. However, this solver may show a very slow convergence, and more efficient techniques include:

\begin{itemize}

\item
{\it Fast solvers}

Since the control problem is in the form of a target problem, a careful discretization preserving {\it causality} would allow for the use of fast marching/fast sweeping solvers. We refer the reader to the discussion and references in \cite{falconeferretilibro}.
 
\item
{\it Modified policy iteration}

In this technique, an inexact policy iteration is implemented by alternating iterations in which a new feedback policy is computed to iterations of plain linear advection. Again, we do not provide details, and rather refer the reader to \cite{FS16} and the references therein, for a detailed study of the hybrid case.

\item
{\it Parallel computing}

In order to apply domain decomposition techniques, the hyperbolic nature of the problem requires some care to optimize the communication between the treads. Recent works on this subject are \cite{CCFP12}, \cite{F16}, \cite{Fh16}. In particular, the \emph{obstacle problem} treated in \cite{Fh16} has a strong link with our case.
\end{itemize}

\section{Numerical tests}\label{Sect:test}

We provide in this section a numerical validation of the technique under consideration, using some typical situations in race strategy. According to the general features of dynamic programming techniques, the optimal control is computed in feedback form: in the semi-Lagrangian scheme \eqref{scheme2}, the optimal continuous control at a certain node is provided by the control $u$ achieving the min, while the need for a switch is indicated by the situation in which
\begin{equation*}
\min \left( NV^\Delta(x_j,q), \Sigma\left(x_j,q,V^\Delta\right)\right) = NV^\Delta(x_j,q),
\end{equation*}
and in this case, by the definition of the transition operator $N$, the optimal switching is towards the dynamics $q'$ achieving the min in \eqref{scheme11}.

\subsubsection*{Test 1}
As a first test we choose the very basic case of a constant average direction of wind ($\bar a=0$), with $(x_1,x_2,x_3)\in[-1.4,1.4]\times[0,2]\times[-1,1]$ taking as target (windward mark) a disc centered in $(0,1.8)$ of radius $0.04$. The goal of this test is to understand the practical effects of the diffusion term $\sigma$ (associated to the stochastic component of the wind direction) on the solution of the navigation problem and give some clues about its interpretation. \\
State constraints have been implemented by penalization of the boundary value, introducing a stopping cost $\bar b=100$, with $\lambda=10^{-6}$, and the tacking cost $C(x,1,2)=C(x,2,1)\equiv 2$. We adopted the simplified dynamics described in Sec. \ref{exampledetailed}, with a constant vessel speed $r\equiv 0.05$. The semi-Lagrangian scheme of Sec. \ref{Sect:scheme} has been used, with uniform discretization steps $\Delta x=0.02$ and $\Delta t=0.1$.

\begin{figure}
\centering
\includegraphics[height=4.8cm]{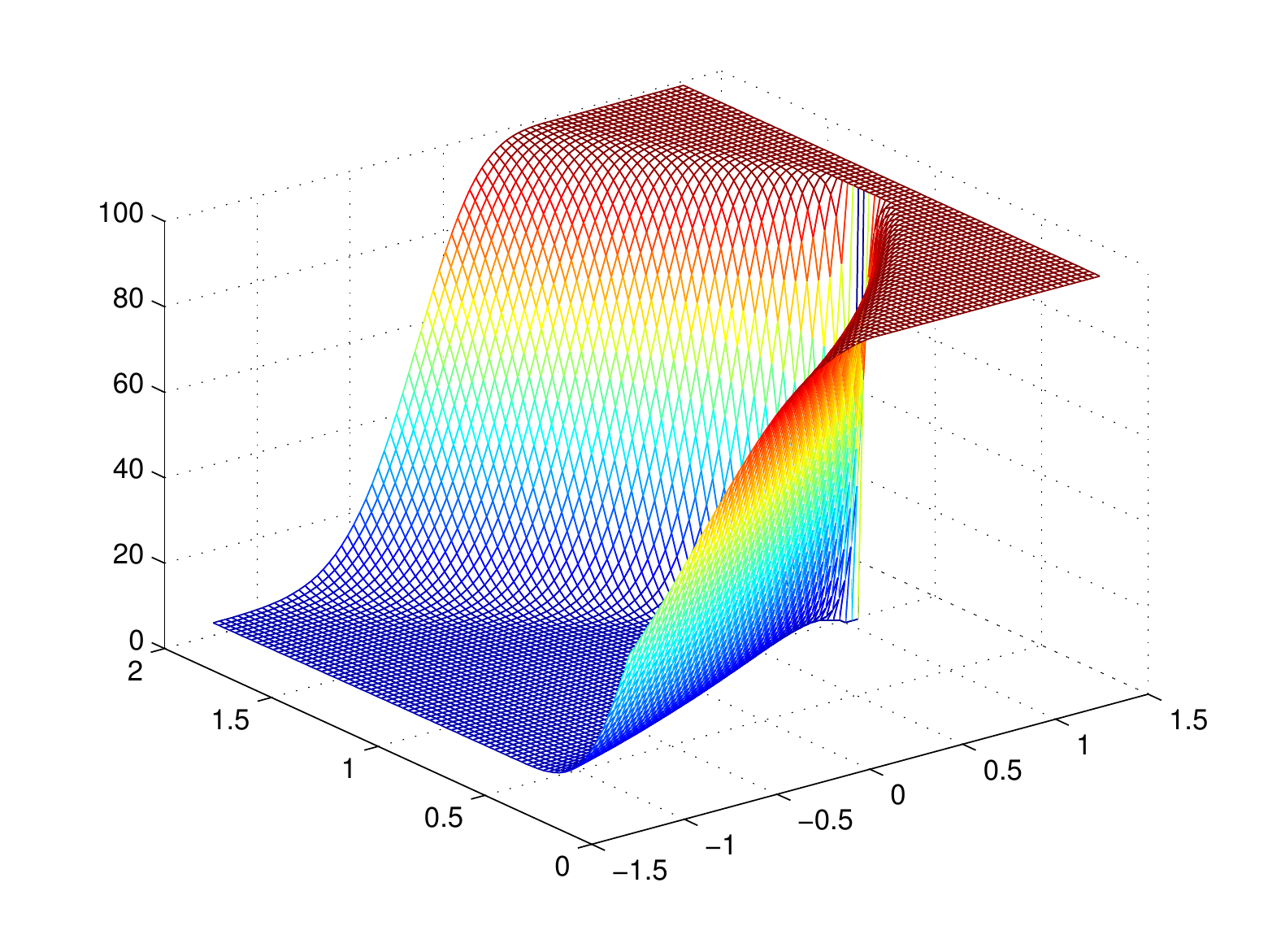}
\includegraphics[height=4.8cm]{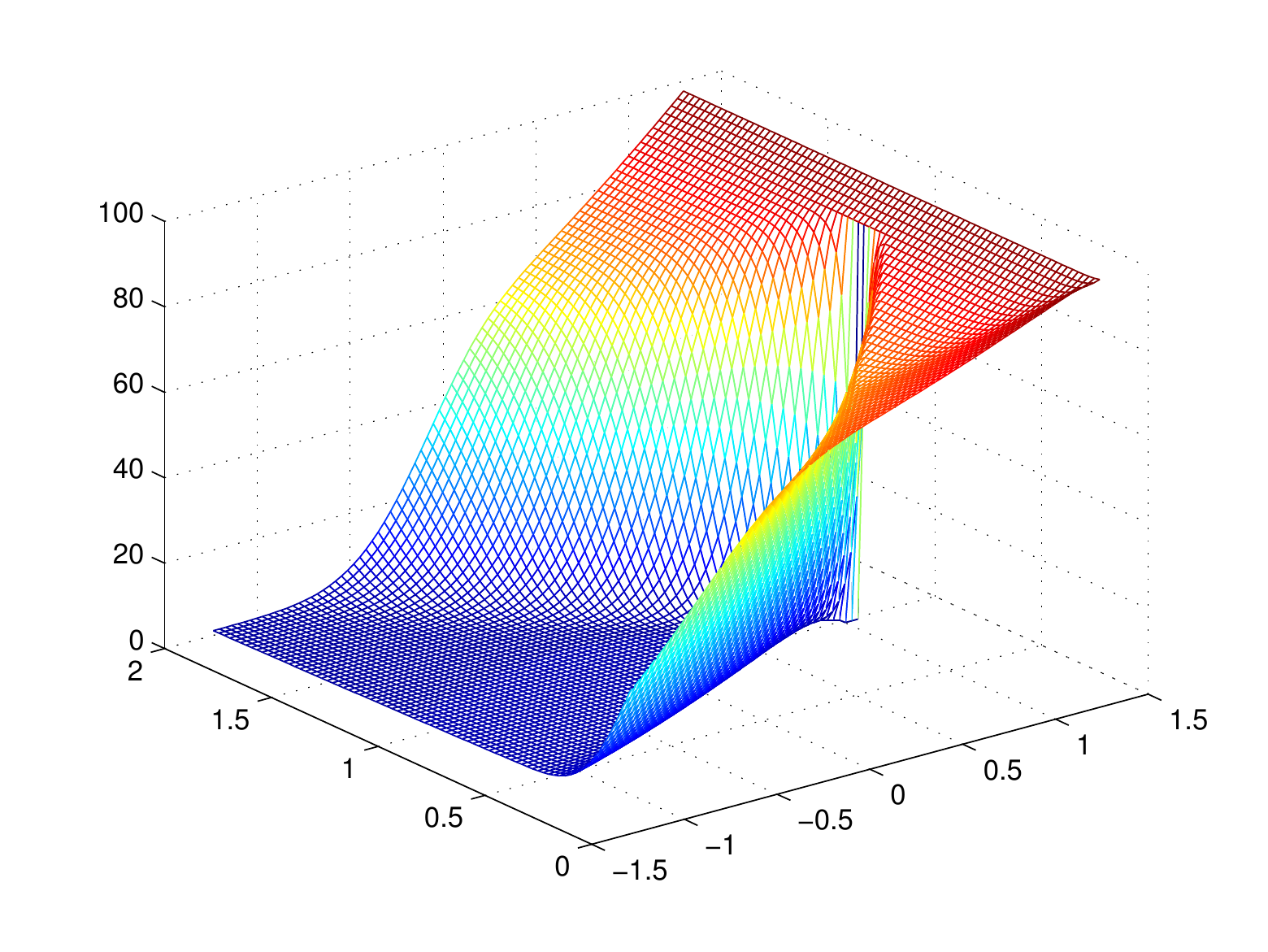} \\
\includegraphics[height=5.2cm]{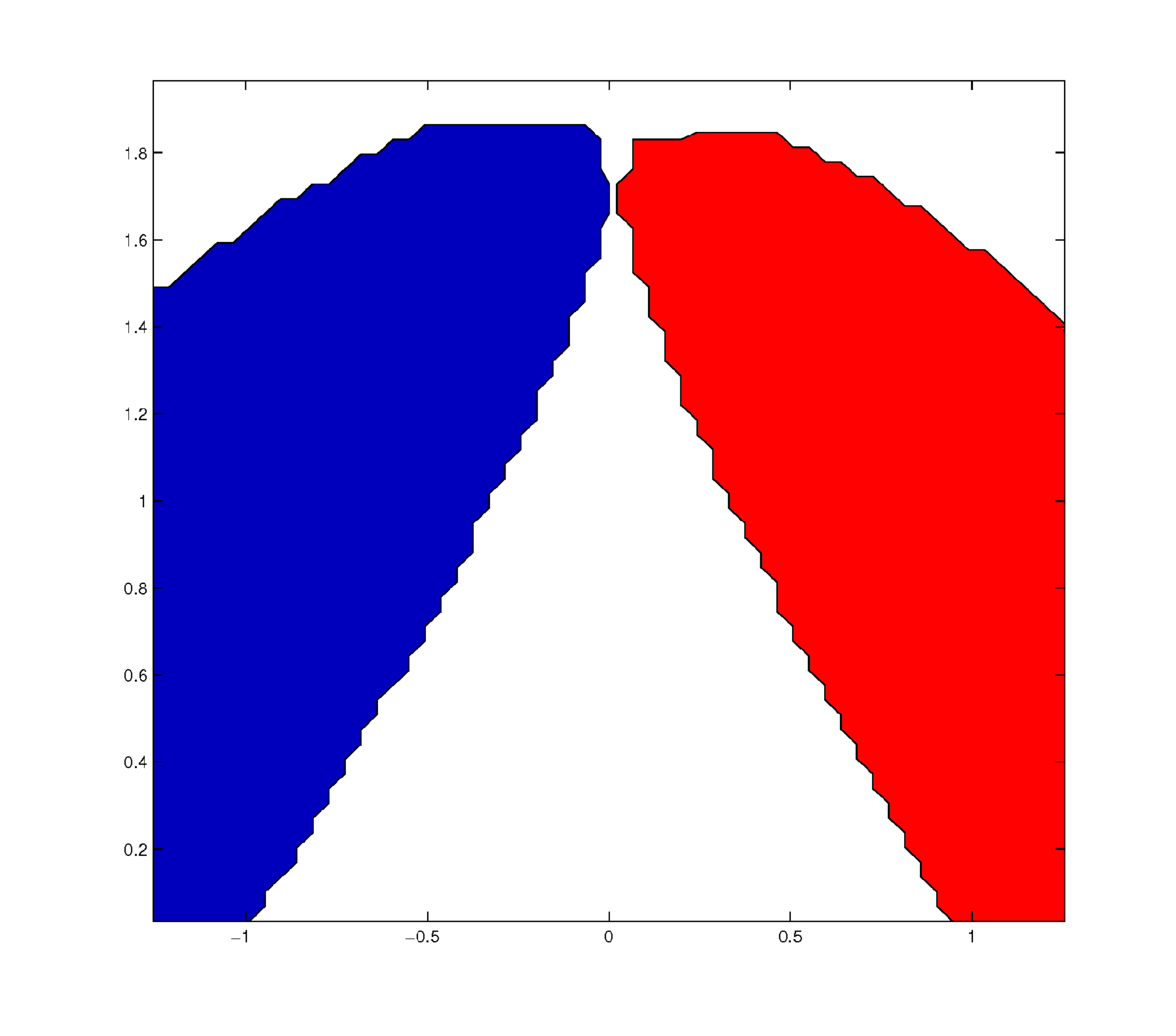}
\includegraphics[height=5.2cm]{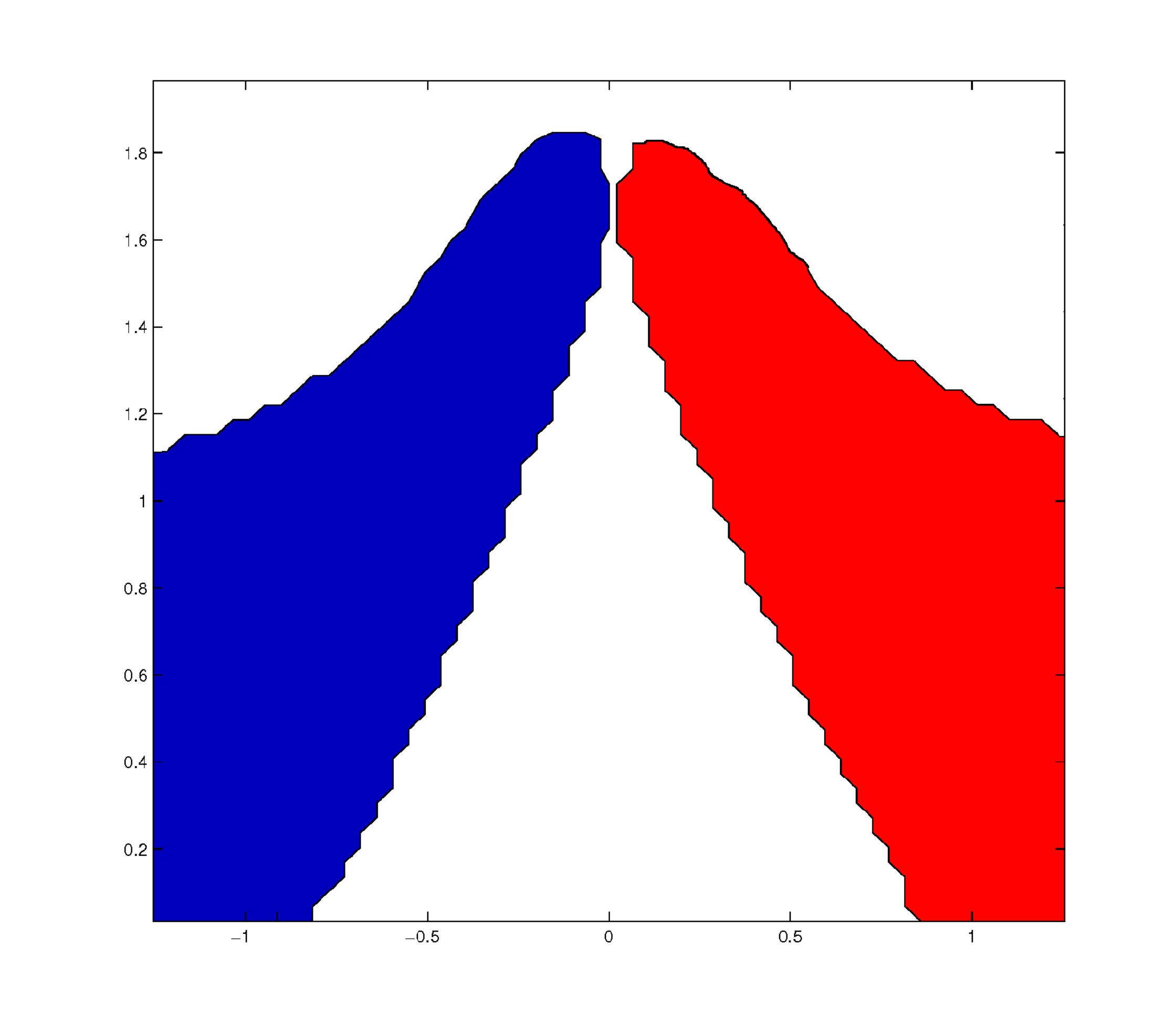}
\caption{Test 1: comparison between the solution of  \eqref{qvi1} in the case of  drift $\bar a=0$, with $\bar \sigma=0$ (upper/left) and $\bar \sigma=0.02$ (upper/right), and corresponding optimal switching maps (lower). Both the comparisons are made for $x_3=0$. The space between the two switching regions is usually called \emph{tacking triangle} (see Fig. \ref{fig:test1theo})} \label{fig:test1plots}
\end{figure}

First, we compare the value functions (associated to the first dynamics \eqref{qvi1}) in the case of $\bar\sigma=0$ and $\bar\sigma=0.02$ (Fig. \eqref{fig:test1plots}). We can observe thet even a small value of $\bar\sigma$ has a strong impact in the regularity of the value function, by smoothing out the sharp gradients of the solution appearing on the so-called {\em lay lines} (i.e., the lines that allow reaching the target without further tacking). The same figure also shows the switching regions for $x_3=0$. Here, the optimal strategy requires tacking to port in the blue region and to starboard in the red region, while no tacking is required in the white region in between. We can observe that a higher stochastic component reduces the width of this region. We will come back shortly to this point.

\begin{figure}
\centering
\includegraphics[height=4.8cm]{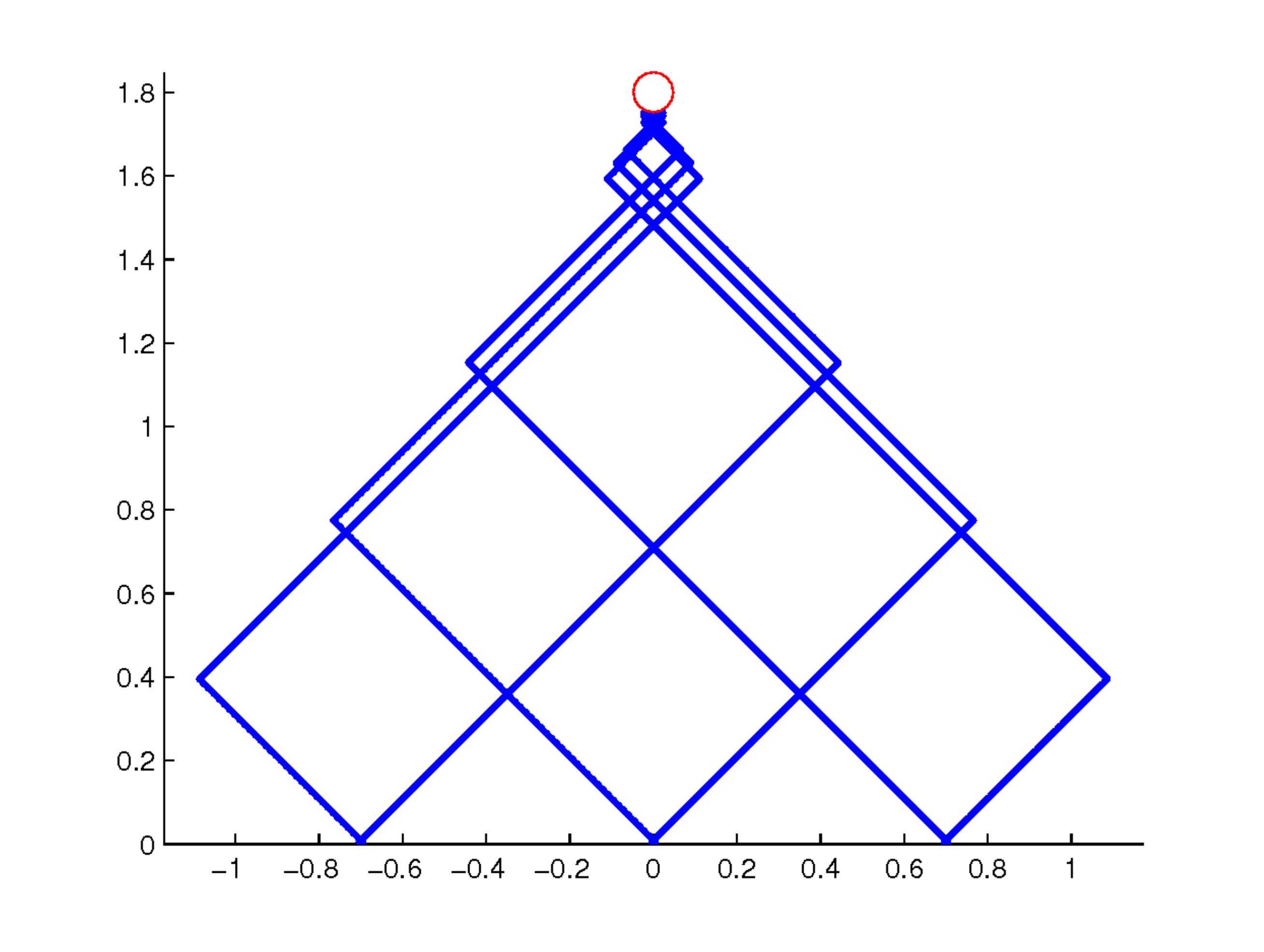}
\includegraphics[height=4.8cm]{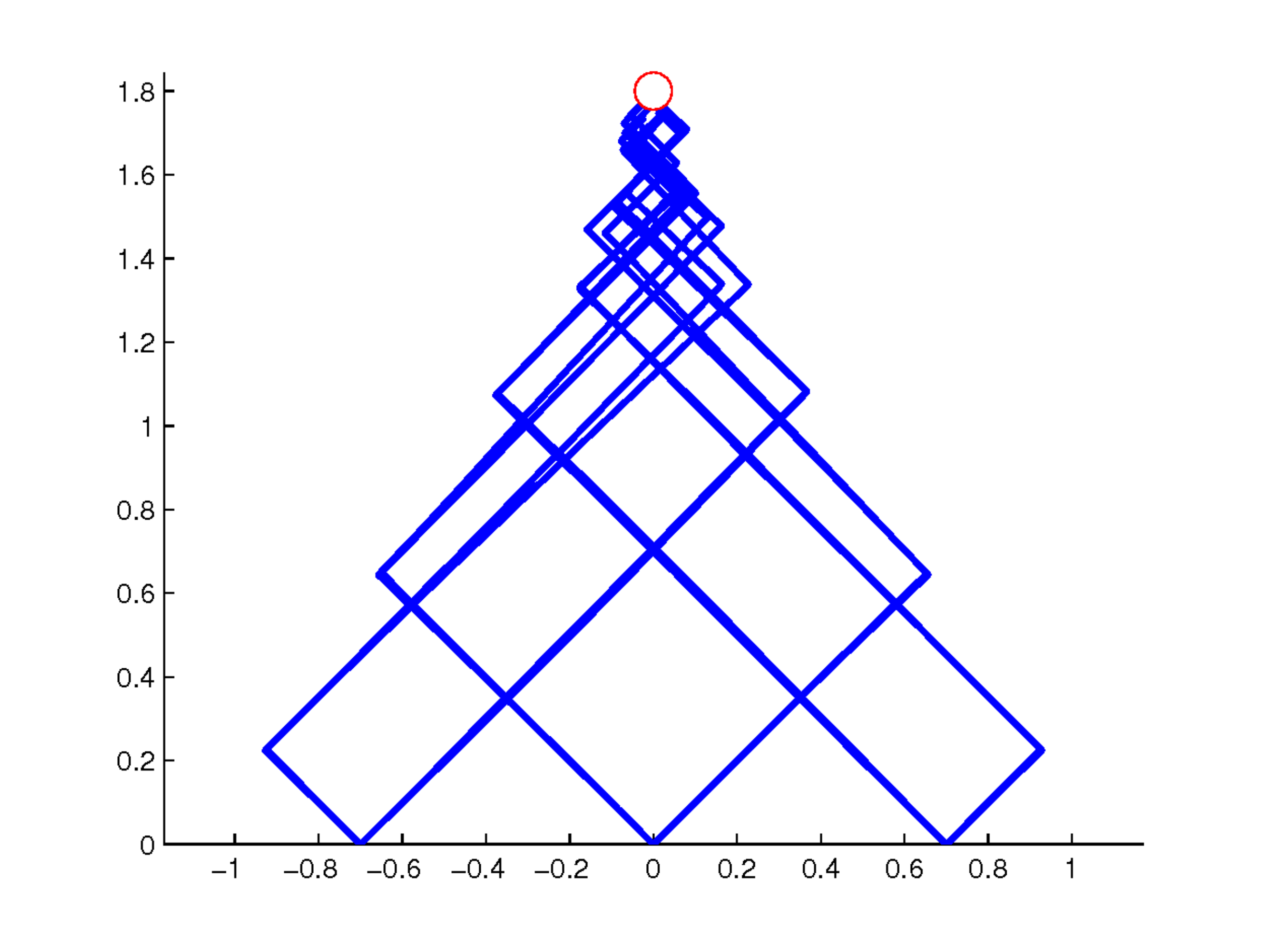} \\
\includegraphics[height=4.8cm]{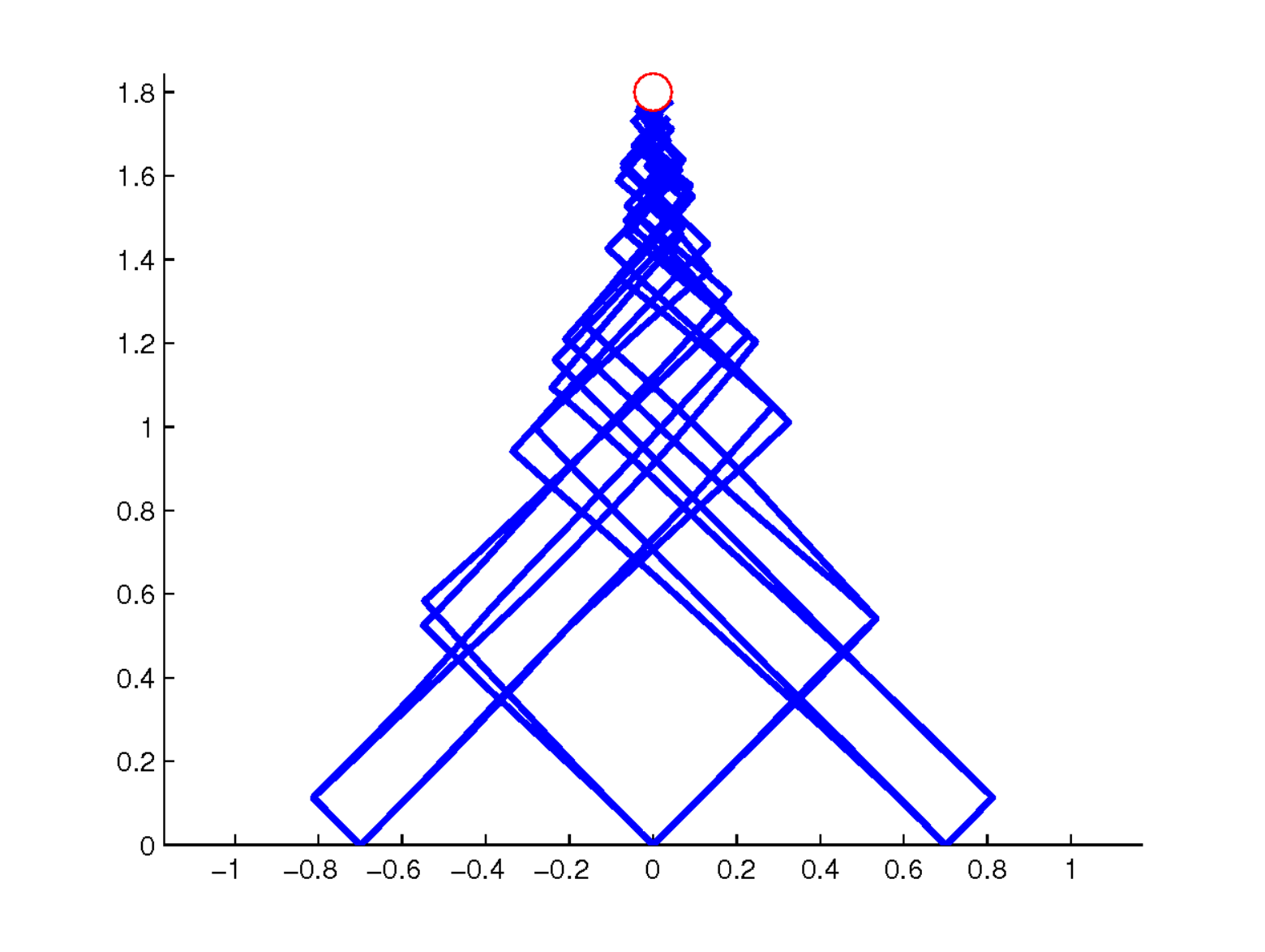}
\includegraphics[height=4.8cm]{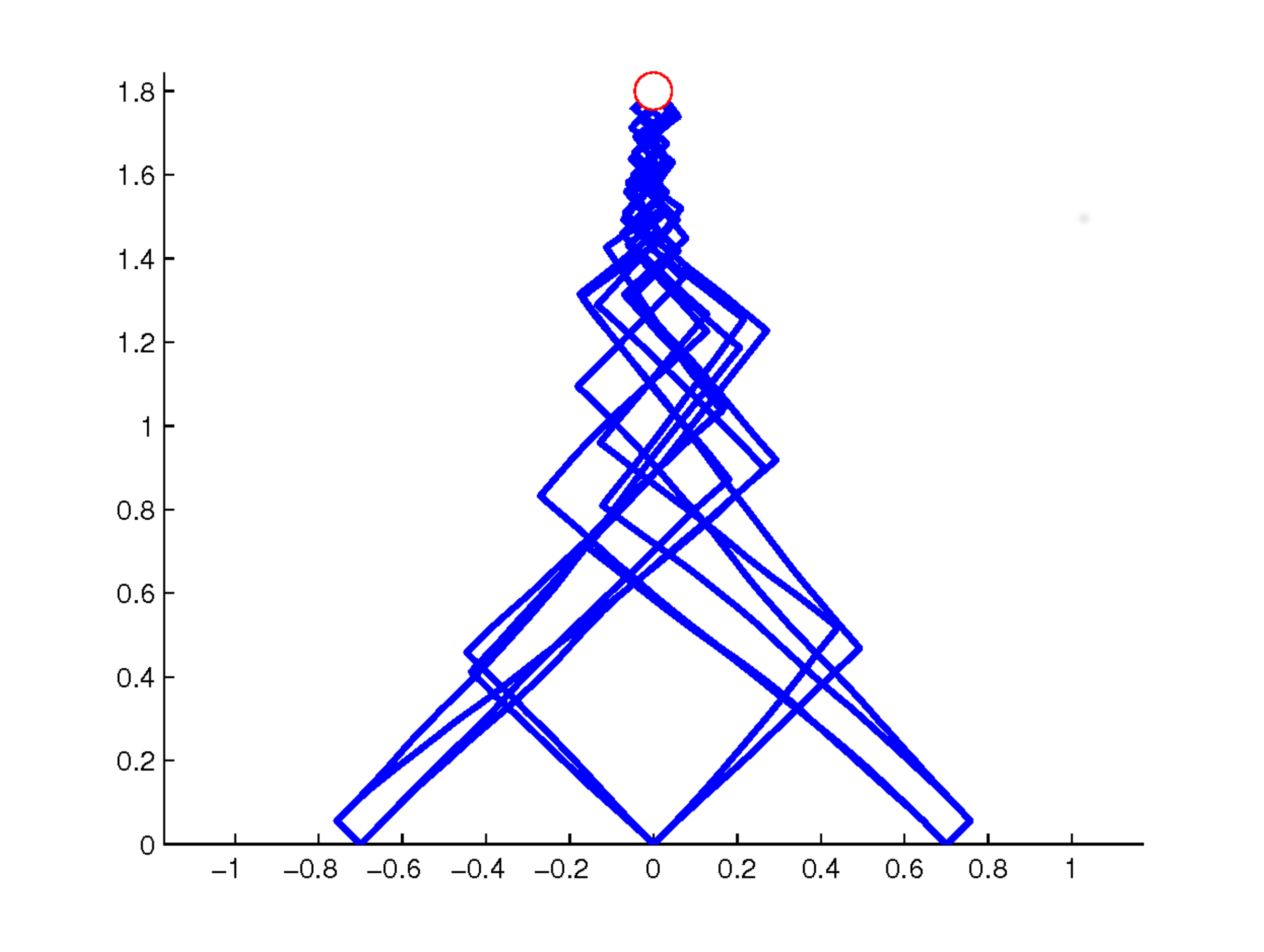}
\caption{Test 1: some sample optimal trajectories for various starting points and different diffusion coefficients. The target is identified by the red circle (above/left $\bar\sigma=0$, above/right $\bar\sigma=0.01$, bottom/left $\bar\sigma=0.05$, bottom/right $\bar\sigma=0.1$)} \label{fig:test1traj}
\end{figure}

In order to simulate the optimal trajectories, we use a standard stochastic Euler scheme (see \cite{KP}) to approximate \eqref{eq:simpl_dyn}, with the increments $\Delta W$ of the Brownian process simulated via a pseudorandom gaussian number generator.
Several tests are performed for different values of $\bar\sigma$. In each case the starting point of the trajectory is set at the points $x=(-0.7,0,0)$, $(0,0,0)$, $(0.7,0,0)$, for $q=1,2$. The results are shown in Fig. \ref{fig:test1plots}. 

In the first case ($\bar\sigma=0$) the trajectories are fully deterministic. In this case, the best strategy is to minimize the number of switches. The strategy pursued is then to wait until the lay lines, switching in this way just once -- in practice, a second switch closer to the target is caused by numerical smoothing effects. The simulation is repeated for values of $\bar\sigma=0.01, 0.05, 0.1$. The number 
of switches in the dynamics increases at the increase of  $\bar\sigma$, i.e., in case of a larger diffusion coefficient for the wind direction, the optimal dynamics prefers to pay more times the cost of switching to remain in the center of the domain. The trajectories tend then to cluster inside a cone that progressively shrinks for increasing $\bar\sigma$.

\begin{figure}[t]
\centering

\includegraphics[height=5.5cm]{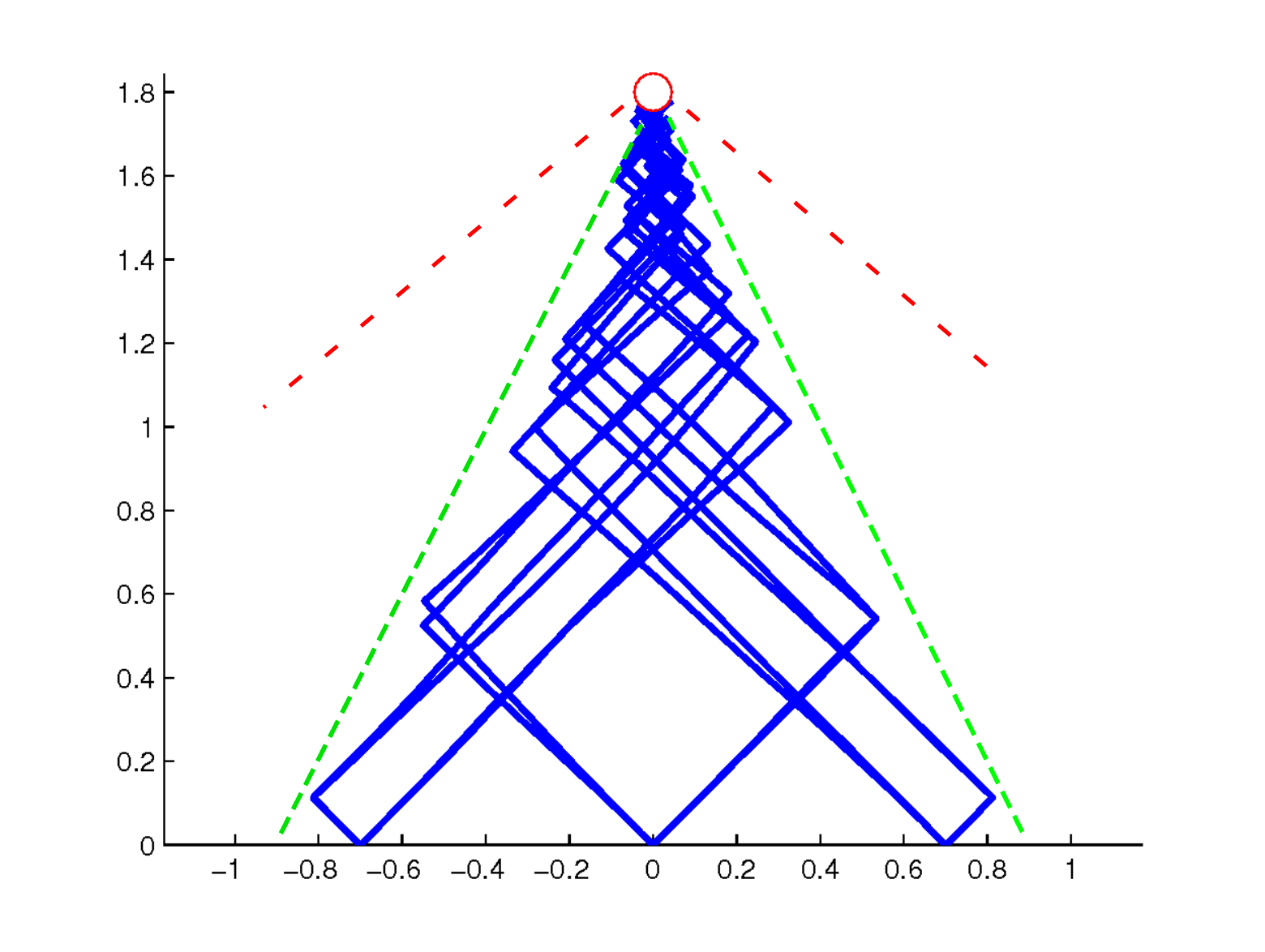}
\includegraphics[height=5.5cm]{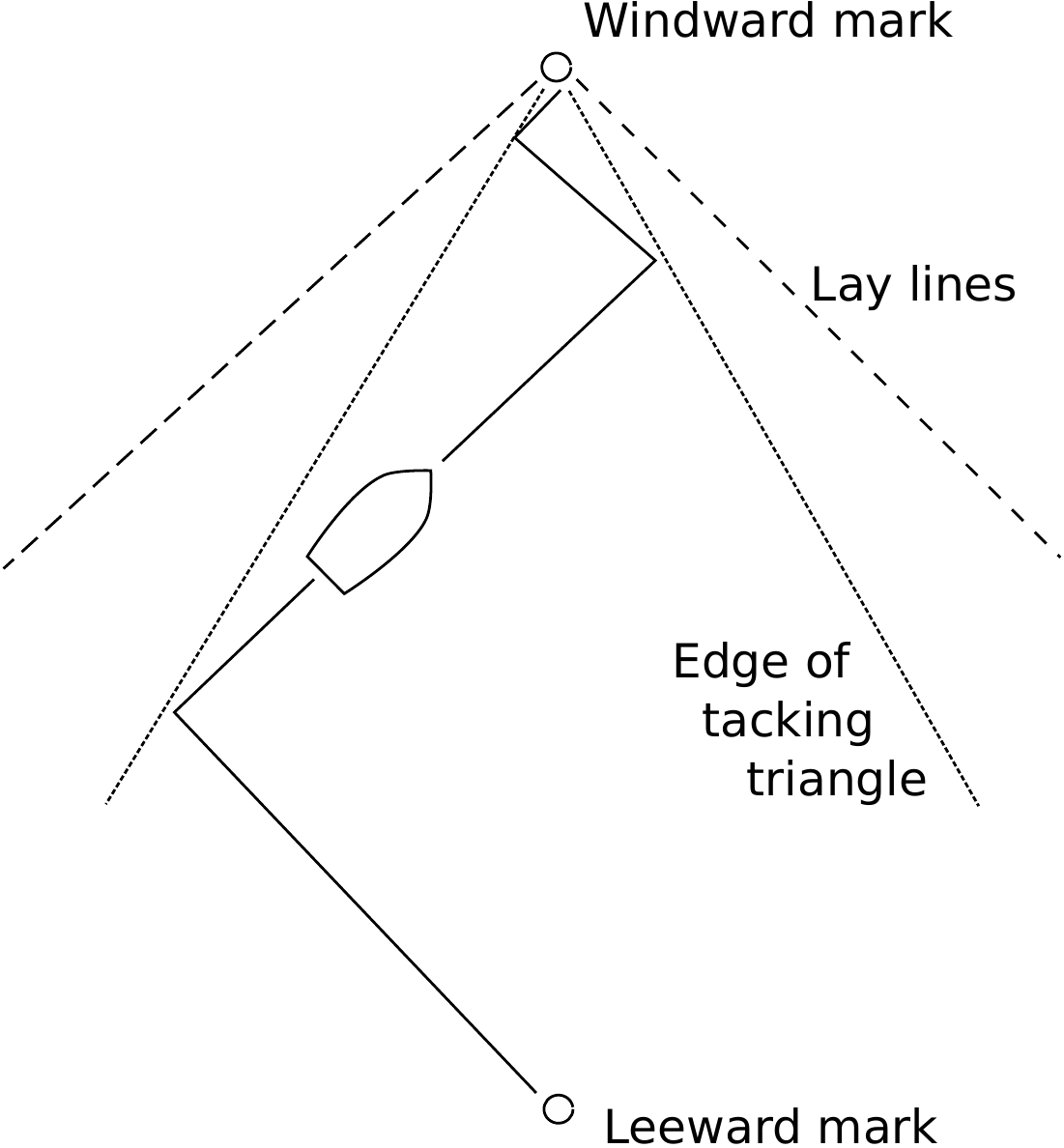}
\caption{Test 1: The results of our tests compared with classic sailing tactics.} \label{fig:test1theo}
\vspace{1cm}
\includegraphics[height=4.7cm]{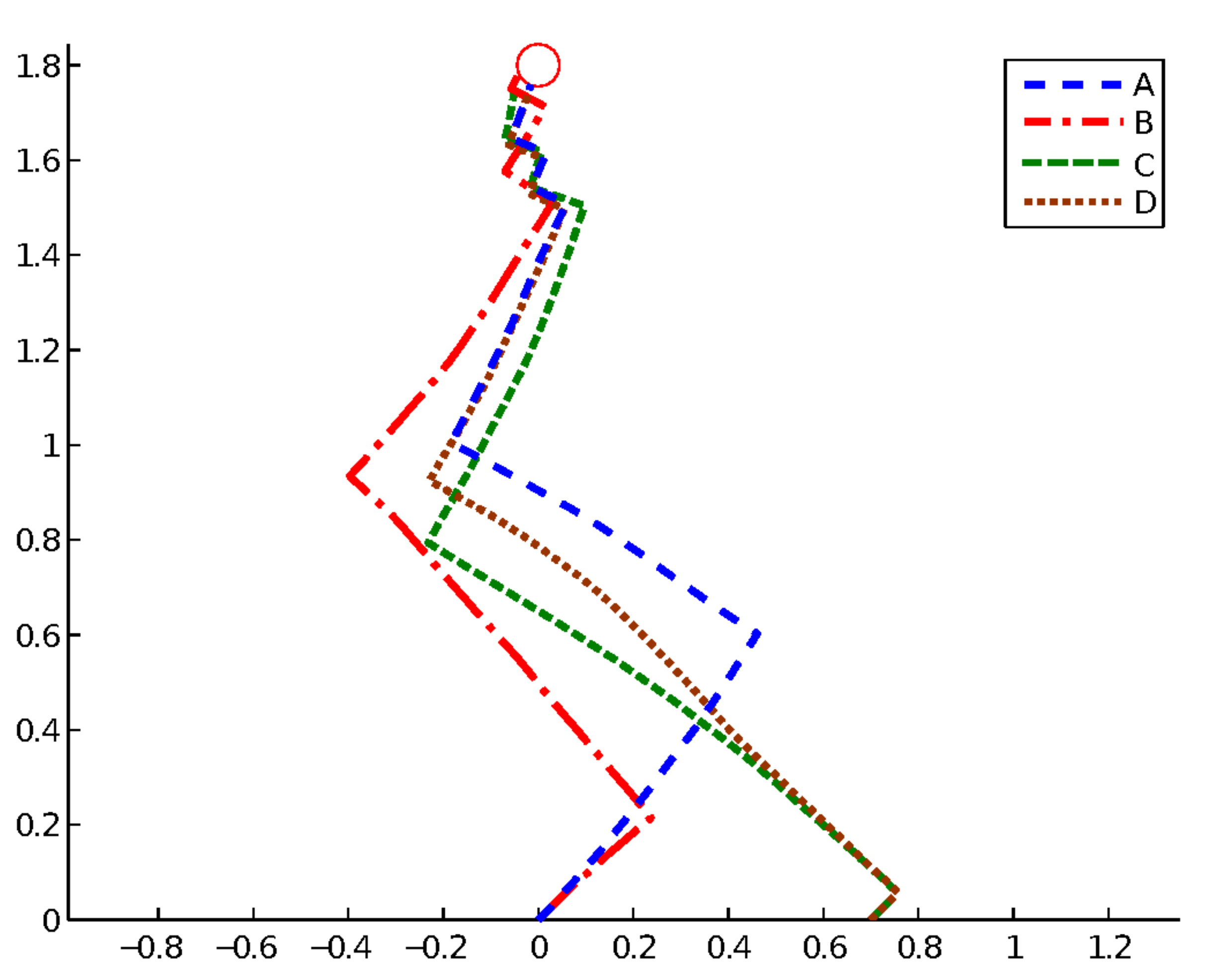}
\hspace{1cm}
\includegraphics[height=4.7cm]{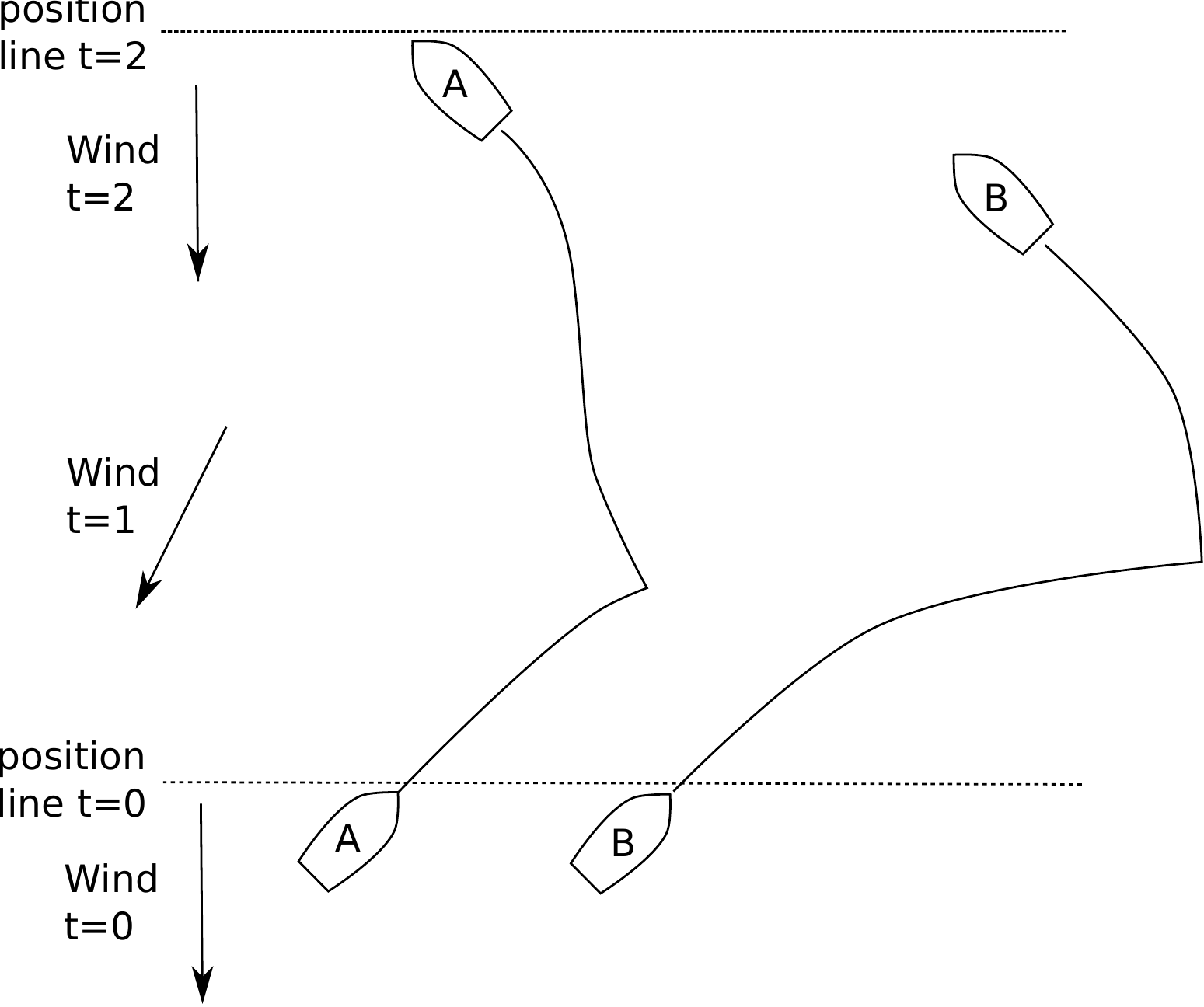}
\caption{Test 1: (left) automatic evaluation of the optimal tacking in presence of wind perturbations and (right) a classical strategy of `tacking on a lift'.} \label{fig:test2comm}
\end{figure}

This effect is well-known in the classical sailing tactics theory \cite{R08}, which suggests (see Fig. \ref{fig:test1theo}) to tack inside a cone, called \emph{tacking triangle}, which coincides with the space between the lay lines in the case of an (ideal) constant wind, and shrinks as the wind becomes more variable. Another qualitative strategy performed in this scenario is known as {\em tacking on a lift} strategy: when the wind direction shows both clockwise and anti-clockwise variations, the former should be preferably used to tack to starboard, while the latter to tack to port. This is schematically shown in Fig. \ref{fig:test2comm} (right), in which A tacks when the wind rotates, and gains on B that does it later.\\
Both tacking strategies are given a more quantitative application with the technique under consideration (cf. Fig. \ref{fig:test2comm}). Without any need to re-compute the value function, the optimal trajectory of the various players is only determined by the particular realization of the wind evolution.

\begin{figure}
\centering
\includegraphics[height=4.8cm]{test1diff005.pdf}
\includegraphics[height=4.8cm]{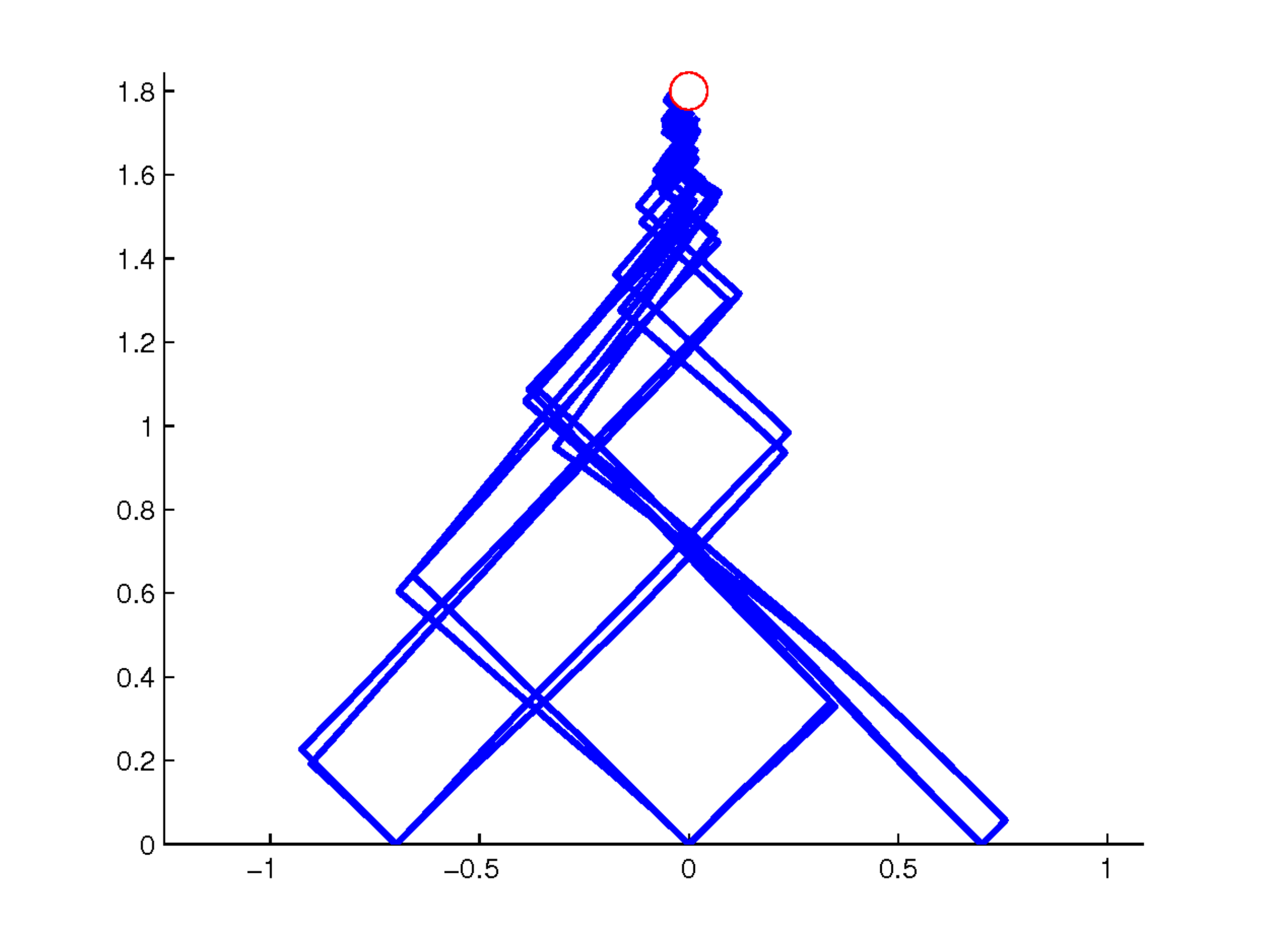} \\
\includegraphics[height=4.8cm]{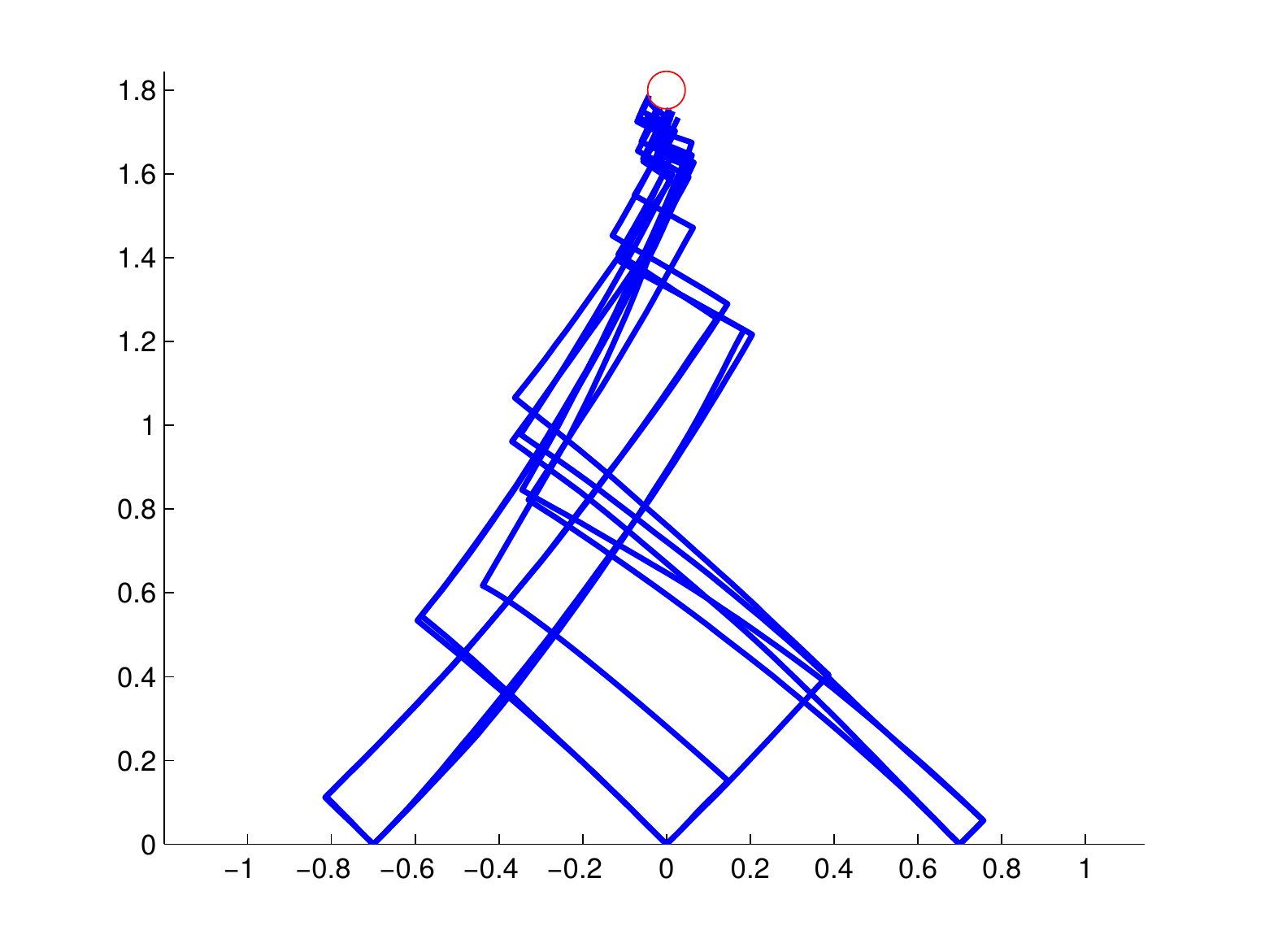}
\includegraphics[height=4.8cm]{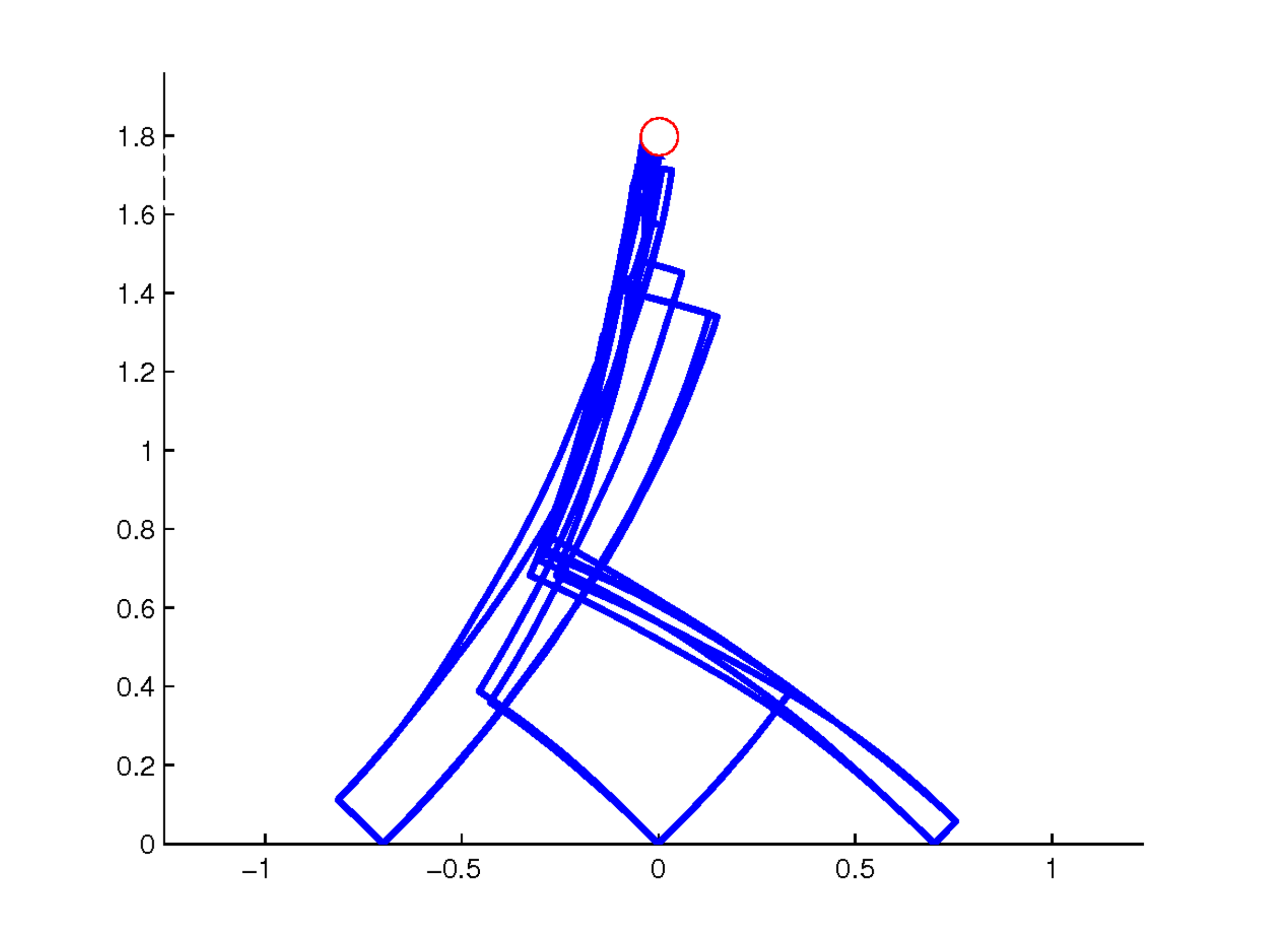}
\caption{Test 2: some sample optimal trajectories for various starting points and different drift values. The diffusion coefficent is fixed $\bar\sigma=0.05$, variable drift (above/left $a=0$, above/right $a=0.05$, bottom/left $a=0.15$, bottom/right $a=0.3$)} \label{fig:test2traj}
\end{figure}

\subsubsection*{Test 2} 
\begin{figure}
\includegraphics[height=4.8cm]{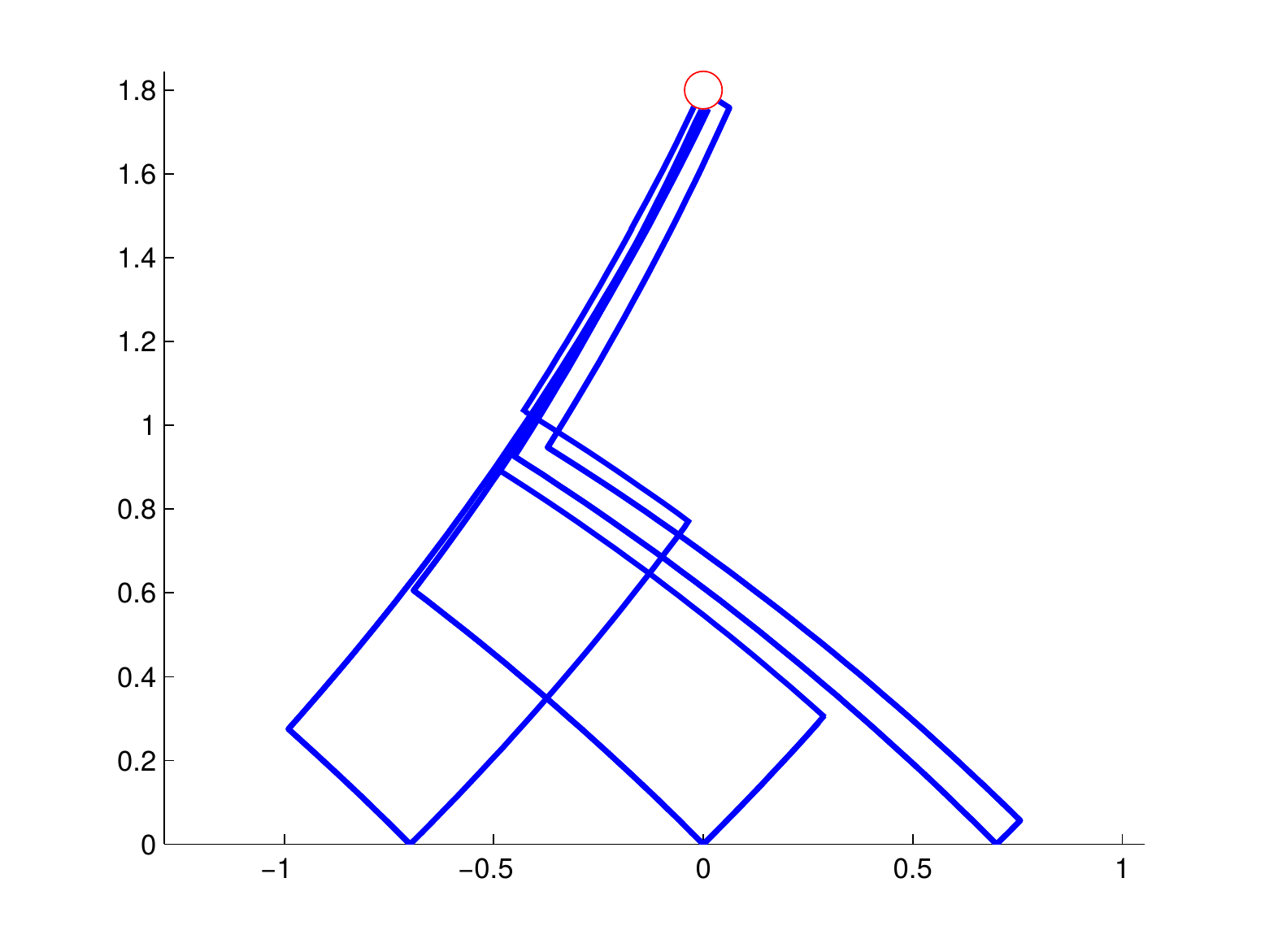}
\includegraphics[height=4.8cm]{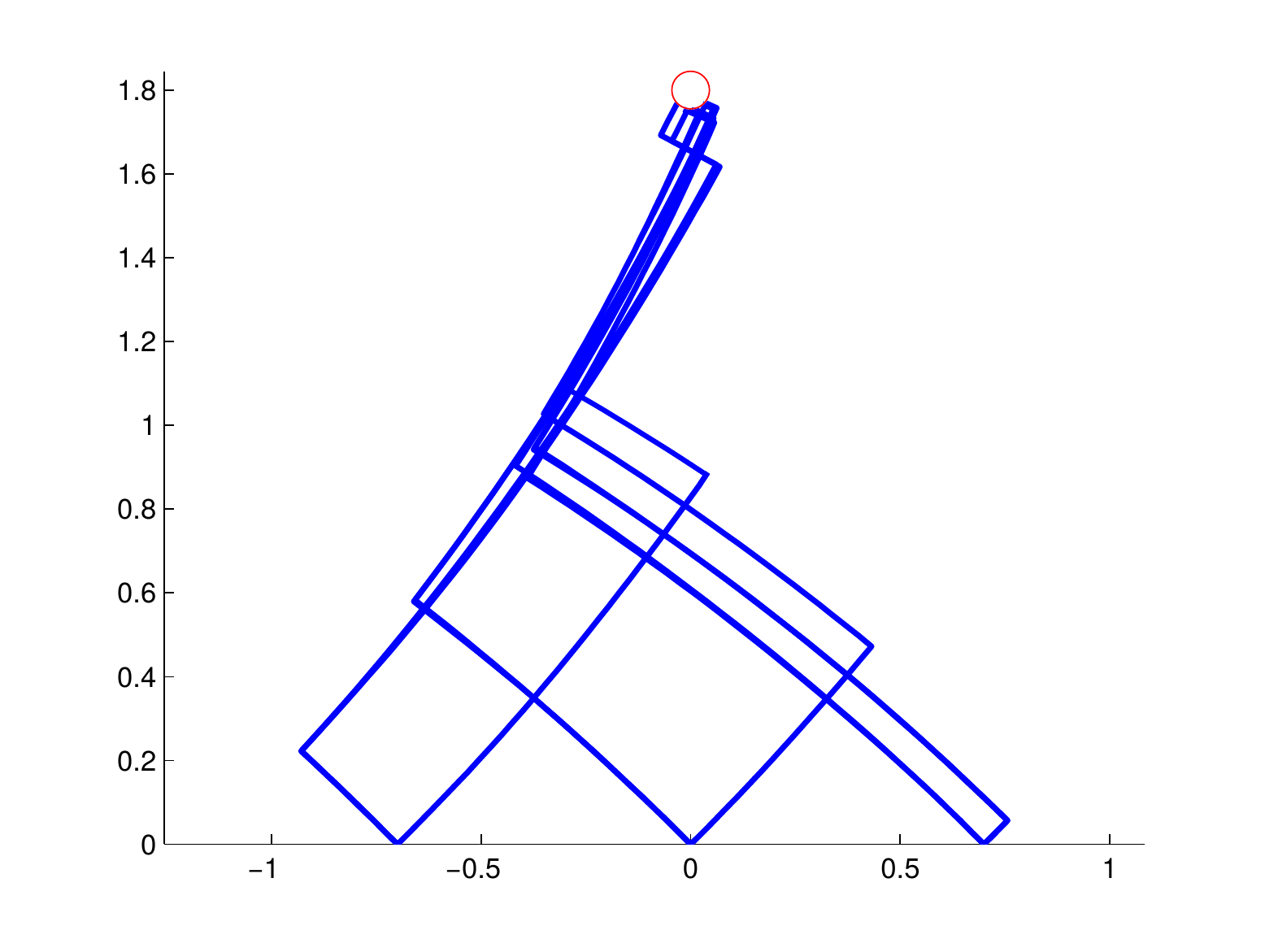} \\
\includegraphics[height=4.8cm]{test2diff005drift015.pdf}
\includegraphics[height=4.8cm]{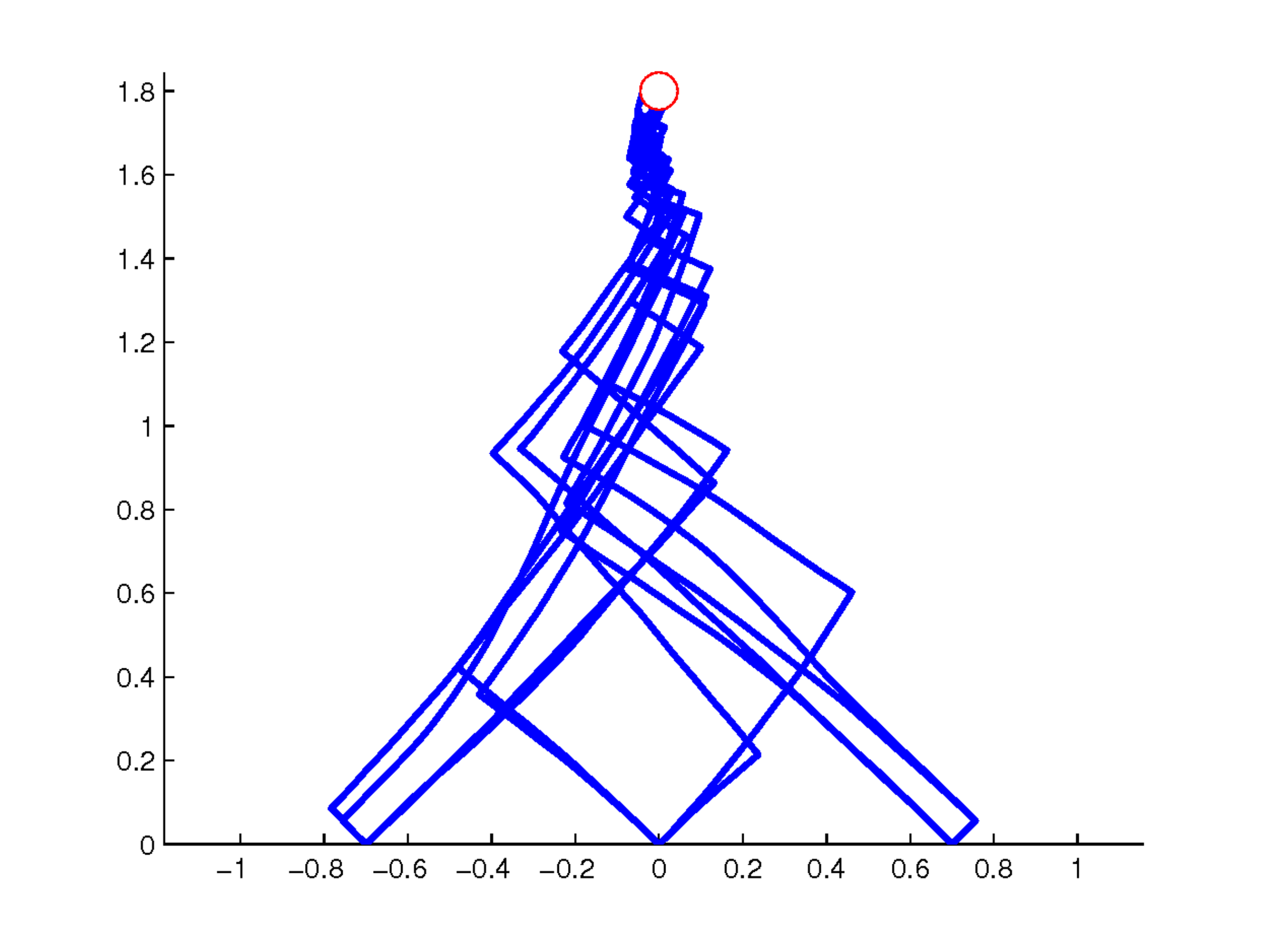}
\caption{Test 2: some sample optimal trajectories for various starting points, constant drift $a=0.15$ and variable diffusion coefficients (above/left $\bar\sigma=0$, above/right $\bar\sigma=0.01$, bottom/left $\bar\sigma=0.05$, bottom/right $\bar\sigma=0.1$)} \label{fig:test2trajdrift}
\end{figure}

In the second test, we impose a nonzero value to the drift $a$ in \eqref{eq:wind} to obtain an average anti-clockwise rotation, possibly changing the diffusion term $\sigma$. In this case, we adopt the complete dynamics of Sec. \ref{exampledetailed} introducing a control variable $u\in [0,\pi/2]$ and a vessel speed given by
\begin{equation}\label{windeq}
r(s,u)\equiv r(u)=0.05\left((\pi/4)^2-(u-\pi/4)^2\right).
\end{equation}
We note that, for this model, the maximum projection of the speed on the axis $x_2$ is attained for $u=\pi/4$, so that, at least qualitatively, the optimal trajectories coincide (see Fig. \ref{fig:test1traj} (left/bottom) and Fig. \ref{fig:test2trajdrift} (left/top)) with the ones obtained with the reduced dynamics. Of course, this is not expected to hold unless the target is purely windward.

Using the same setting of Test 1, we first compare the optimal trajectories obtained for various values of the drift $\bar a$, for a constant value of diffusion $\bar\sigma=0.05$. As in the previous test, we consider trajectories starting from the points $(-0.7,0,0)$, $(0,0,0)$, $(0.7,0,0)$ and $q=1,2$. The results are shown in Fig. \ref{fig:test2traj}. In the second set of tests, we repeat the same scenario fixing a positive drift $\bar a=0.15$ with diffusion coefficients $\bar \sigma=0,0.01, 0.05, 0.1$. The results are reported in Fig. \ref{fig:test2trajdrift}. 

Observe that, as $a>0$, the approximate symmetry of optimal trajectory in the previous test is progressively lost. In practice, if the average rotation of the wind is anti-clockwise, then the best strategy is to occupy the left region of the domain (i.e., for $x_1$ small or negative), and this behavior is enhanced by large values of $a$. This global strategy include also the `tacking on a lift' strategy, which holds also in the situation of a variable average wind direction with stochastic variations.

The computed optimal solutions correctly blends the two strategies. The tendency to occupy the left side of the state domain is higher for higher values of $a$ and lower values of $\sigma$.

\subsubsection*{Test 3} In this test we show how the model can handle state constraints that naturally appear in applications (presence of obstacles, coasts, etc). Among the wide literature on constrained feedback control, we quote the recent works \cite{altarovici2013general, aubin2008dirichlet} and the references therein. In our case, the presence of obstacles is taken into account with a penalization on the speed of motion. Then, considered an obstacle $\Gamma\subset \R^d$, we replace the function $r(s,u)$ in \eqref{eq:motion} with the following state dependent velocity
\begin{equation*}
r(x,s,u):=
\begin{cases}
r(s,u)  & \text{ if }x\in \R^d\setminus \Gamma,\\
0 & \text{ if }x\in \Gamma.
\end{cases}
\end{equation*}
with $r(s,u)$ as in \eqref{windeq}.

\begin{figure}
\centering
\includegraphics[height=5cm]{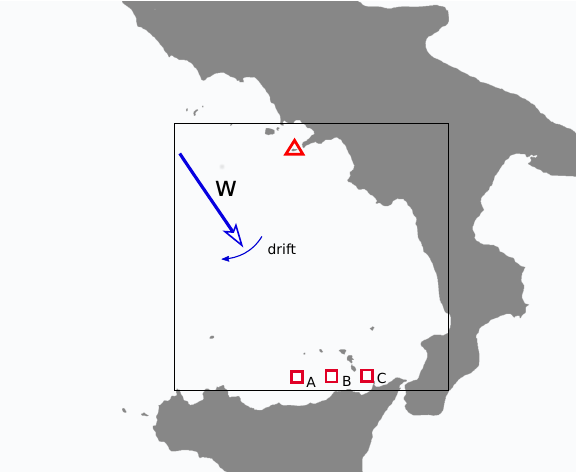}
\includegraphics[height=5cm]{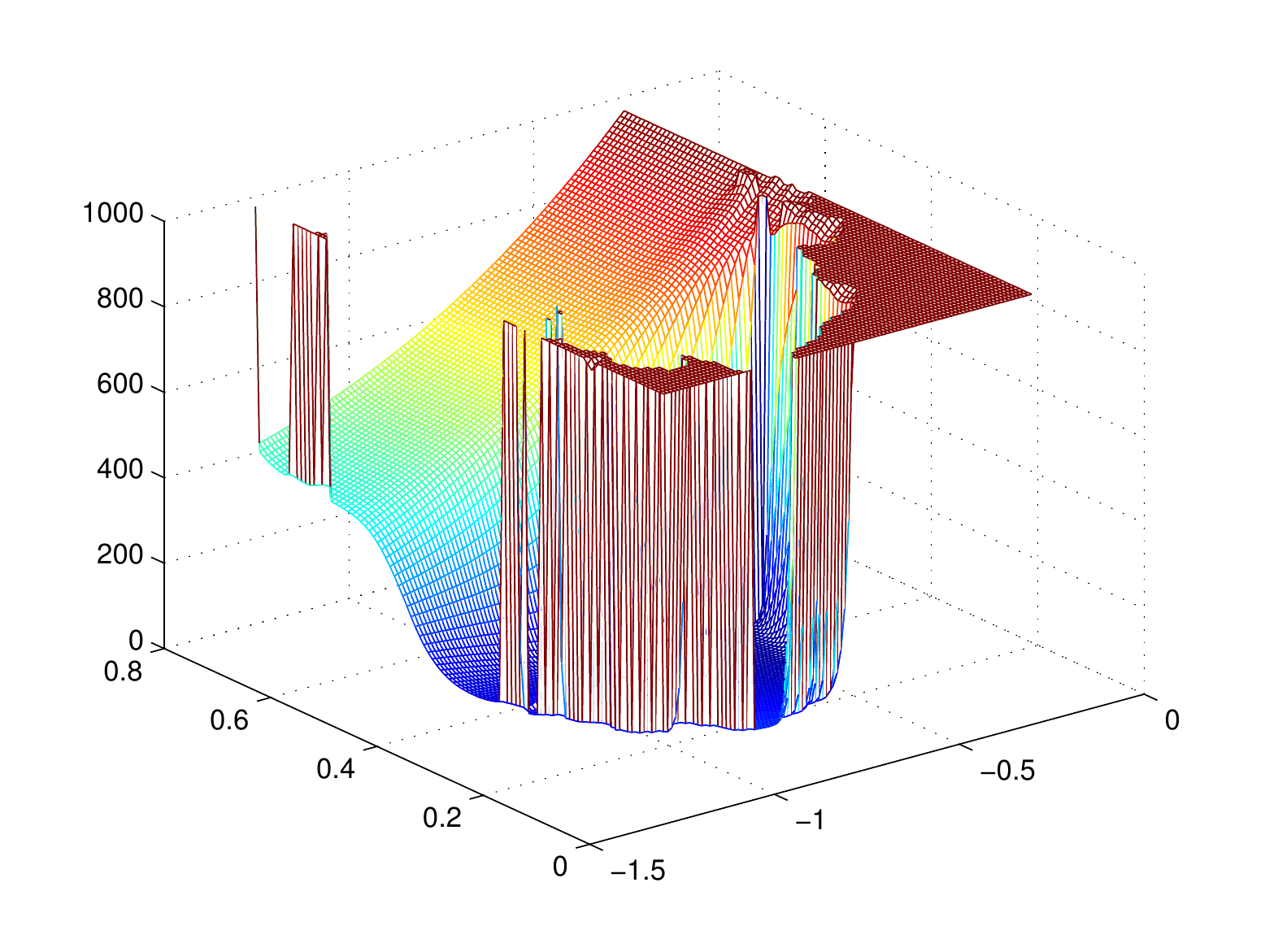} 
\caption{Test 3: (left) the area of interest (black rectangle), the objective (red triangle) and the starting points (red squares). The initial wind direction (blue arrow marked with a 'W') is shown as the direction of rotation of the drift (small blue arrow marked 'drift'). (Right) The value function (dynamics $q=1$) for $x_3=-0.5$, $\bar \sigma=0.005$. } \label{fig:test3}
\bigskip

\includegraphics[height=5cm]{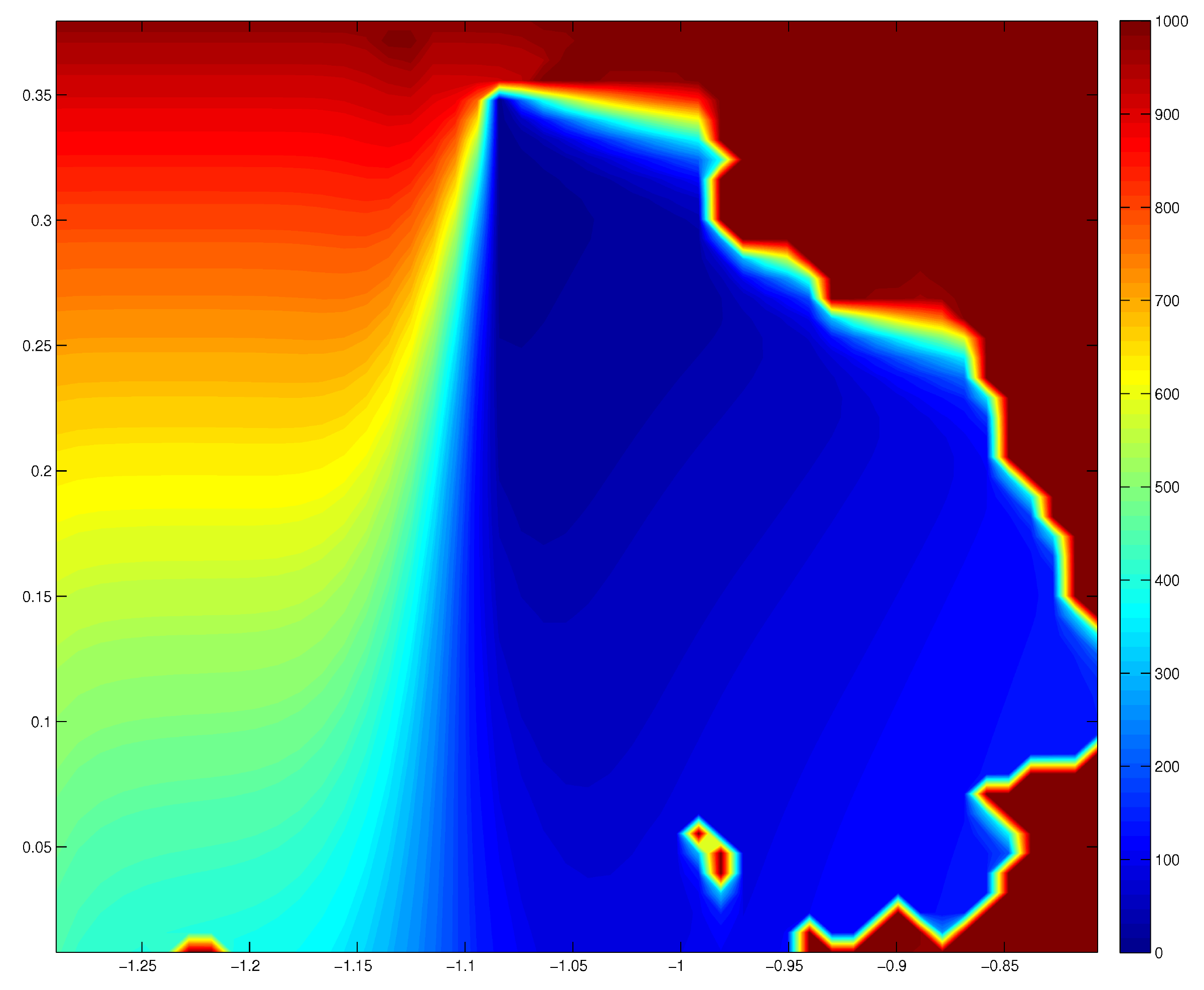}
\includegraphics[height=5cm]{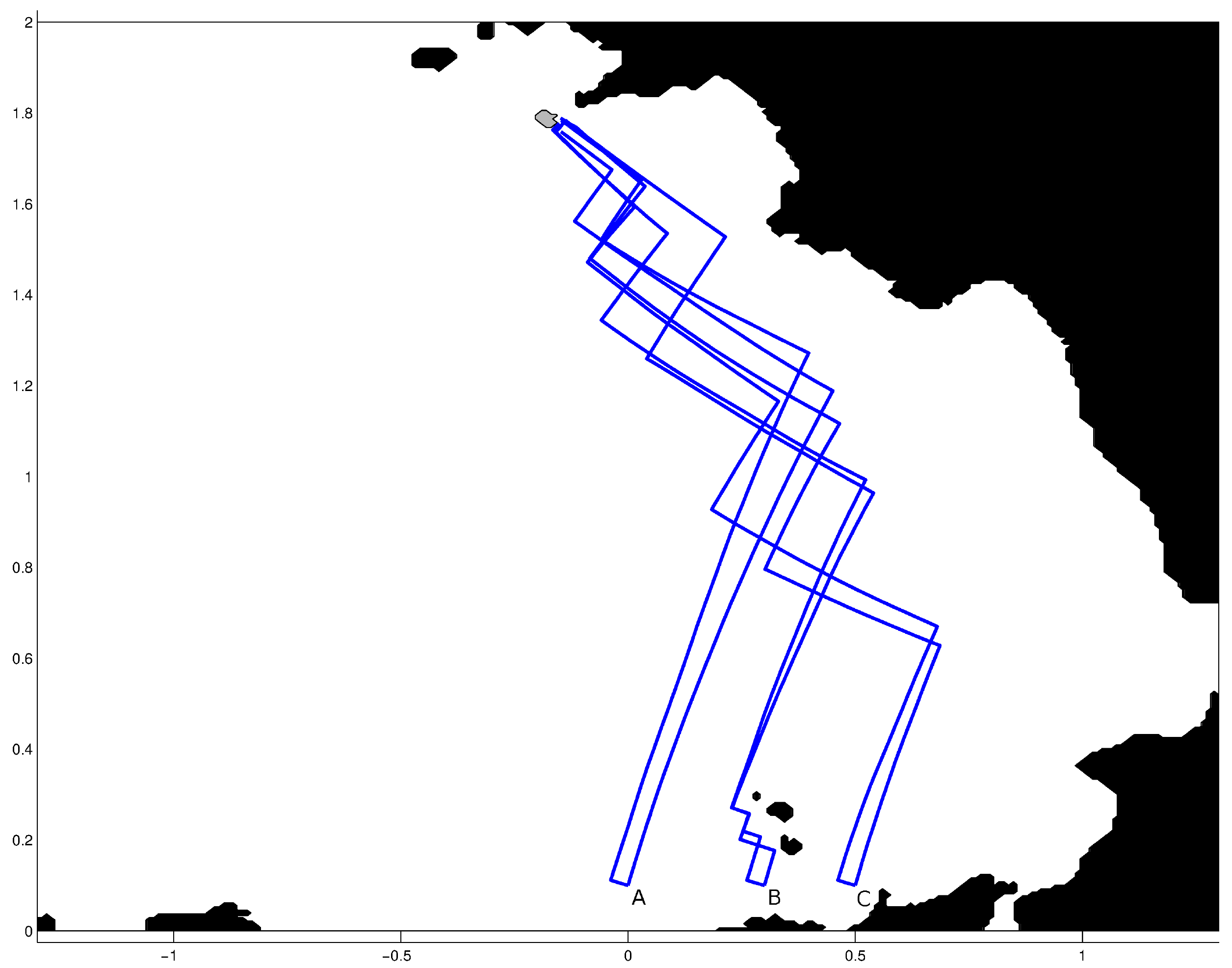} 
\caption{Test 3: (left) contour plot of the value function (dynamics $q=1$) for $x_3=-0.5$, $\sigma=0.005$ and (right) sample optimal trajectories starting from the points of interest.} \label{fig:test3traj}
\end{figure}

In this example we still keep a simplified wind evolution with a wind field with initial direction $\theta(0)=0.5$, speed $\bar s=1$, drift $a=-0.15$ and diffusion coefficient $\bar \sigma=0.05$, while the cost parameters are chosen as in the previous tests. In Fig. \ref{fig:test3} we show the computational domain, which represents the southern area of Tyrrhenian sea\footnote{Map by Free Vector Maps {https://freevectormaps.com/italy/IT-EPS-01-0001} }. We set target at the island of Capri (red triangle in Fig. \ref{fig:test3}), and consider trajectories starting from the points A, B, C, (marked with red squares in Fig. \ref{fig:test3}) close to the coast of Sicily. Note that a change in the starting point does not require to recompute the optimal solution, which has already been obtained in feedback form.

The negative drift models an average clockwise rotation of the wind, which involve a preference for the right side of the domain. On the other hand, this strategy is limited by the presence of the coast. 

Fig. \ref{fig:test3traj} (left) shows a contour plot of the value function of the problem for $x_3=-0.5$ and $q=1$, where the constraint $\Gamma$ is apparent. In Fig. \ref{fig:test3traj} (right) we show some optimal trajectories starting from the points A, B, C. Here, the presence of the coast causes a less extreme strategy than in the case of Test 2. Note that, in the trajectory obtained starting from the point B, some small islands obstruct the way to the `natural' choice, and must be avoided by the player before pointing rightwards.

\section{Acknowledgments}
This work was partially supported by the New Frontiers Grant NST 0001 by the Austrian Academy of Sciences, the Haute-Normandie Regional Council via the M2NUM project, by the INdAM--GNCS project ``Metodi numerici per equazioni iperboliche e cinetiche e applicazioni'' and by Roma Tre University.

\end{document}